\documentclass[a4paper,11pt]{article}
\usepackage[latin1]{inputenc}
\usepackage[T1]{fontenc}
\usepackage[francais]{babel}
\usepackage{amssymb}
\usepackage{amsfonts, amsmath}
\newtheorem{defn}{D\'efinition}[section]

\newtheorem{prop}[defn]{Proposition}

\newtheorem{cor}[defn]{Corollaire}
\newtheorem{rem}{Remarque}[section]

\numberwithin{equation}{section}
\title{FONCTIONS ET INTEGRALES ELLIPTIQUES}
\author{\textbf{A. Lesfari}
\\\emph{Département de Mathématiques}
\\\emph{Faculté des Sciences}
\\\emph{Université Chouaïb Doukkali}
\\\emph{B.P. 20, El-Jadida, Maroc}.
\\\emph{E. mail : lesfariahmed@yahoo.fr, lesfari@ucd.ac.ma}}
\date{}

\begin{document}
\maketitle
\begin{abstract}
This paper presents the basic ideas and properties of elliptic
functions and elliptic integrals as an expository essay. It
explores some of their numerous consequences and includes
applications to some problems such as the simple pendulum, the
Euler rigid body motion and some others integrable hamiltonian systems.
\end{abstract}
\vskip0.6cm
\textbf{Plan :}\\
1. Fonctions elliptiques.\\
2. Fonctions de Weierstrass.\\
\null\hskip0.6cm2.1. Fonction $\wp$ de Weierstrass.\\
\null\hskip0.6cm2.2. Fonction $\zeta$ de Weierstrass.\\
\null\hskip0.6cm2.3. Fonction $\sigma$ de Weierstrass.\\
3. Intégrales elliptiques et fonctions de Jacobi.\\
4. Applications.\\
\null\hskip0.6cm4.1. Le pendule simple.\\
\null\hskip0.6cm4.2. Le corps solide d'Euler.\\
\null\hskip0.6cm4.3. Une famille de systèmes intégrables.\\
\null\hskip0.6cm4.4. Equations aux dérivées partielles
non-linéaires de
Schrödinger.\\
\null\hskip0.6cm4.5. Le champ de Yang-Mills avec groupe de jauge $SU(2)$.\\
\null\hskip0.6cm4.6. Appendice.\\

\vskip0.6cm Les fonctions elliptiques interviennent dans des
domaines très divers. Le but de ce travail est de montrer quelques
résultats fondamentaux sur ces fonctions et de les appliquer à des
situations concrètes.

\section{Fonctions elliptiques}

Soient $\omega_1$ et $\omega_2$ deux nombres complexes,
$\mathbb{R}$-linéairement indépendants (c'est-à-dire tels que
$\omega_2$ n'est pas nulle et que le quotient
$\frac{\omega_1}{\omega_2}$ ne soit pas réel ou ce qui revient au
même que la partie imaginaire $\mbox{Im}
\frac{\omega_1}{\omega_2}$ du rapport $\frac{\omega_1}{\omega_2}$
n'est pas nulle). On considère le réseau
\begin{eqnarray}
\Lambda&=&\mathbb{Z}\omega_1\oplus\mathbb{Z}\omega_2,\nonumber\\
&=&\{\omega\equiv m\omega_1+n\omega_2 : m,n\in
\mathbb{Z}\},\nonumber
\end{eqnarray}
c'est un sous-groupe discret de $\mathbb{C}$ et il forme un
ensemble de parallélogrammes.
\begin{defn}
On appelle parallélogramme fondamental engendré par $\omega_1$ et
$\omega_2$ tout parallélogramme $\Pi$ de sommets d'affixes $z_0$,
$z_0+\alpha\omega_1$, $z_0+\alpha\omega_2$,
$z_0+\alpha\omega_1+\beta\omega_2$ avec $z_0\in\mathbb{C}$, $0\leq
\alpha, \beta\leq 1$. Autrement dit, il est défini par le compact
$$\Pi=\{z_0+\alpha\omega_1+\beta\omega_2\ : z_0\in\mathbb{C},\quad \alpha,
\beta\in [0,1]\}.$$
\end{defn}
Le quotient de $\mathbb{C}$ par la relation d'équivalence
déterminée par $\Lambda$ :
$$z_1, z_2\in \mathbb{C},\quad z_1\sim z_2
 \mbox{mod}. \Lambda \Longleftrightarrow z_1-z_2\in \Lambda,$$ est
un tore noté $\mathbb{C}/\Lambda$. Celui-ci est homéomorphe à
$S^1\times S^1$, visualisable par le recollement deux à deux des
côtés d'un carré ou parallélogramme.

\begin{defn}
On dit qu'une fonction $f$ de $\mathbb{C}$ dans $\mathbb{C}$ est
doublement périodique de périodes $\omega_1$ et $\omega_2$, si et
seulemnt si,
\begin{eqnarray}
f(z+\omega_1)&=&f(z),\nonumber\\
f(z+\omega_2)&=&f(z).\nonumber
\end{eqnarray}
Autrement dit, si et seulement si,
$$f(z+\omega)=f(z),\quad \forall \omega \in \Lambda.$$
On dit aussi que $f$ est $\Lambda$-périodique.
\end{defn}

\begin{rem}
Les éléments $\omega_1$ et $\omega_2$ ne sont pas uniques. Plus
précisément, si $\omega_1$ et $\omega_2$ sont deux périodes de
$f,$ alors $-\omega_1$ et $-\omega_2$ sont également deux périodes
de $f$ et toute période de $f$ s'écrit sous la forme
$$\omega=m\omega_1+n\omega_2,\quad m,n\in \mathbb{Z}.$$
En effet, pour les entiers positifs c'est évident. Pour les
entiers négatifs, on a pour $k=1, 2$ :
\begin{eqnarray}
f(z-\omega_k)&=&f((z-\omega_k)+\omega_k),\nonumber\\
&=&f(z).\nonumber
\end{eqnarray}
\end{rem}

\begin{defn}
On dit qu'une fonction $f$ de $\mathbb{C}$ dans $\mathbb{C}$ est
elliptique si et seulemnt si elle est méromorphe et doublement
périodique.
\end{defn}

\begin{prop}
Il n'existe pas de fonction elliptique $f$ non constante qui soit
holomorphe sur $\mathbb{C}$. Autrement dit, toute fonction
elliptique $f$ n'ayant pas de pôles est une constante.
\end{prop}
\emph{Démonstration}: Si $f$ n'a pas de pôles, alors elle est
bornée dans le parallélogramme fondamental $\Pi$ car celui-ci est
compact. Or la fonction $f$ est doublement périodique, donc elle
est bornée sur $\mathbb{C}$ car tout point de $\mathbb{C}$ se
ramène à un point de $\Pi$ en lui appliquant une translation du
réseau. Par conséquent, $f$ est constante en vertu du théorème de
Liouville. $\square$

\begin{rem}
On déduit de la proposition précédente qu'une fonction elliptique
non constante possède au moins un pôle dans le parallélogramme
fondamental.
\end{rem}

\begin{prop}
Toute fonction elliptique non constante a un nombre fini de pôles
et un nombre fini de zéros dans un parallélogramme fondamental.
\end{prop}
\emph{Démonstration}: Rappelons qu'un point
$a\in\Omega\subset\mathbb{C}$ est un point d'accumulation s'il
existe une suite $(z_n)_{n\in\mathbb{N}}$ d'éléments de
$D\setminus\{a\}$ telle que : $\lim_{n\rightarrow\infty}z_n=a$.
Soit $P(f)\equiv f^{-1}(\{\infty\})$ l'ensemble des pôles de la
fonction $$f
:D\longrightarrow\overline{\mathbb{C}}=\mathbb{C\cup\{\infty\}}.$$
Comme $f$ est méromorphe, alors l'ensemble $P(f)$ n'admet pas de
point d'accumulation. Donc $f$ a un nombre fini de pôles car sinon
$P(f)$ doit contenir le point limite (point d'accumulation) et
celà est impossible car le point d'accumulation des pôles est une
singularité essentielle. Soit maintenant $$Z(f)=\{b\in D :
f(b)=0\},$$ l'ensemble des zéros de $f$. Comme $f$ n'est pas
constante, l'ensemble $Z(f)$ n'admet pas de point d'accumulation.
Dès lors, pour tout point $b\in Z(f)$, il existe un unique entier
positif $N$ tels que : $$f(z)=(z-b)^Ng(z),$$ où $g$ est une
fonction holomorphe sur $D$ avec $g(b)\neq 0$. En fait l'ensemble
$Z(f)$ est au plus dénombrable. $\square$

\begin{rem}
Comme remarque, notons que le nombre de zéros et de pôles d'une
fonction elliptique non constante est dénombrable. En effet,
l'ensemble des parallélogrammes fondamentals forme un recouvrement
dénombrable de $\mathbb{C}$ et le résultat découle de la
proposition précédente.
\end{rem}

\begin{prop}
Soit $f$ une fonction elliptique et désignons par $b_1,...,b_m$
les pôles de $f$ (chaque pôle étant compté avec multiplicité),
alors
$$\sum_{k=1}^m\mbox{Rés}(f,b_{k})=0.$$
\end{prop}
\emph{Démonstration}: Soit
$\gamma=\gamma_1\cup\gamma_2\cup\gamma_3\cup\gamma_4$ la frontière
du parallélogramme fondamental $\Pi$ relativement au réseau
$\Lambda$, avec
\begin{eqnarray}
\gamma_1&=&[z_0,z_0+\omega_1],\nonumber\\
\gamma_2&=&[z_0+\omega_1,z_0+\omega_1+\omega_2],\nonumber\\
\gamma_3&=&[z_0+\omega_1+\omega_2,z_0+\omega_2],\nonumber\\
\gamma_4&=&[z_0+\omega_2,z_0].\nonumber
\end{eqnarray}
Supposons tout d'abord que $f$ n'a pas de pôles sur la frontière
$\gamma$. D'après le théorème des résidus, on a
\begin{eqnarray}
\sum_{k=1}^m\mbox{Rés}(f,b_{k})&=&\frac{1}{2\pi i}\int_\gamma
f(z)dz,\nonumber\\
&=&\frac{1}{2\pi i}(\int_{\gamma_1} f(z)dz+\int_{\gamma_2}
f(z)dz+\int_{\gamma_3} f(z)dz+\int_{\gamma_4} f(z)dz).\nonumber
\end{eqnarray}
En vertu de la périodicité de $f$ et des sens opposés de
l'intégrale de $f$ sur $\gamma_1$ et $\gamma_2$, on a
\begin{eqnarray}
\int_{\gamma_3}f(z)dz&=&\int_{[z_0+\omega_1+\omega_2,z_0+\omega_2]} f(z)dz,\nonumber\\
&=&\int_{[z_0+\omega_1,z_0]} f(u+\omega_2)du,\quad u\equiv z-\omega_2,\nonumber\\
&=&\int_{[z_0+\omega_1,z_0]} f(u)du,\nonumber\\
&=&-\int_{[z_0,z_0+\omega_1]} f(u)du,\nonumber\\
&=&-\int_{\gamma_1} f(z)dz.\nonumber
\end{eqnarray}
De même, on a
$$\int_{\gamma_4}f(z)dz=-\int_{\gamma_2}f(z)dz,$$
et par conséquent $$\sum_{k=1}^m\mbox{Rés}(f,b_{k})=0.$$ Passons
maintenat au cas où il y'a des pôles sur la frontière $\gamma$ du
parallélogramme fondamental $\Pi$. Alors dans ce cas, on considère
un autre parallélogramme proche de $\Pi$ contenant tous les pôles
se trouvant dans $\Pi$ et de telle façon que sa frontière ne
contienne plus de pôles. On peut toujours, d'après la proposition
1.5, obtenir ce parallélogramme (et donc sa frontière) par
translation du sommet d'affixe $z_0$ de $\Pi$. Le reste consiste à
utiliser un raisonnement similaire au précédent. $\square$

\begin{rem}
Notons que d'après la remarque 1.3 et la proposition précédente,
il n'existe pas de fonction elliptique de premier ordre, i.e., une
fonction elliptique ne peut pas avoir un pôle simple dans un
parallélogramme fondamental. Elle doit avoir au moins deux pôles
simples ou au moins un pôle non simple dans un parallélogramme
fondamental. En effet, avec les notations de la proposition
précédente, si $m=1$ alors celà signifie que la fonction $f$ a un
pôle simple dans le parallélogramme fondamental, ce qui contredit
le résultat de la proposition.
\end{rem}

\begin{rem}
L'ensemble des fonctions elliptiques par rapport à $\Lambda$ est
un sous corps du corps des fonctions méromorphes (la somme, le
produit et le quotient de deux fonctions elliptiques de mêmes
périodes est une fonction elliptique). En dérivant l'expression
$$f(z+\omega)=f(z),\quad \forall \omega \in \Lambda,$$
on obtient
$$f^{(n)}(z+\omega)=f^{(n)}(z),\quad \forall \omega \in \Lambda,$$
ce qui montre que la dérivée $n^{\mbox{ième}}$ d'une fonction
elliptique est aussi une fonction elliptique.
\end{rem}

\begin{prop}
Soit $f$ une fonction elliptique non constante. Désignons par
$a_1,...,a_l$ les zéros de $f$ de multiplicité $n_1,...,n_l$
respectivement et par $b_1,...,b_m$ les pôles de $f$ de
multiplicité $p_1,...,p_m$ respectivement. Alors
$$\sum_{k=1}^ln_k=\sum_{k=1}^m p_k.$$
Autrement dit, le nombre de zéros d'une fonction elliptique non
constante est égal au nombre de ses pôles dans le parallélogramme
fondamental.
\end{prop}
\emph{Démonstration}: D'après le principe de l'argument, on a
\begin{eqnarray}
\mbox{Nombre de zéros de} f-\mbox{Nombre de pôles de}
f&=&\frac{1}{2\pi
i}\int_\gamma \frac{f'(z)}{f(z)}dz,\nonumber\\
&=&\sum_{k=1}^m\mbox{Rés}(\frac{f'}{f},b_{k}),\nonumber\\
&=&0,\nonumber
\end{eqnarray}
en vertu de la proposition 1.6 car d'après la proposition
précédente $\frac{f'(z)}{f(z)}$ est une fonction elliptique et a
les mêmes périodes que $f(z)$. Par conséquent, $\sum_{k=1}^l
n_k=\sum_{k=1}^m p_k$. $\square$

\begin{prop}
Soit $f$ une fonction elliptique. Désignons par $a_1,...,a_l$ les
zéros de $f$ de multiplicité $n_1,...,n_l$ respectivement et par
$b_1,...,b_m$ les pôles de $f$ de multiplicité $p_1,...,p_m$
respectivement. Alors
$$\sum_{k=1}^l n_ka_k-\sum_{k=1}^mp_kb_k=\mbox{période}.$$
\end{prop}
\emph{Démonstration}: Rappelons que si une fonction $\varphi(z)$
est holomorphe dans un domaine $D\subset \mathbb{C}$ et continue
sur $\overline{D}$, alors
$$\frac{1}{2\pi i}\int_\gamma
\varphi(z)\frac{f'(z)}{f(z)}dz=\sum_{k=1}^ln_k\varphi(a_k)-\sum_{k=1}^np_k\varphi(b_k).$$
On pose dans la suite $\varphi(z)=z$ et on utilise les mêmes
notations et arguments de la preuve de la proposition 1.7. Donc
\begin{eqnarray}
\sum_{k=1}^ln_ka_k-\sum_{k=1}^np_kb_k&=&\frac{1}{2\pi
i}\int_\gamma z\frac{f'(z)}{f(z)}dz,\nonumber\\
&=&\frac{1}{2\pi i}\sum_{j=1}^4\int_{\gamma_j}
z\frac{f'(z)}{f(z)}dz.\nonumber
\end{eqnarray}
On a
\begin{eqnarray}
\int_{\gamma_1} z\frac{f'(z)}{f(z)}dz+\int_{\gamma_3}
\zeta\frac{f'(\zeta)}{f(\zeta)}d\zeta&=&\int_{\gamma_1}
z\frac{f'(z)}{f(z)}dz-\int_{\gamma_1}
\zeta\frac{f'(\zeta)}{f(\zeta)}d\zeta,\quad \zeta\equiv z+\omega_2,\nonumber\\
&=&\int_{\gamma_1} (z-\zeta)\frac{f'(z)}{f(z)}dz,\nonumber\\
&=&-\omega_2\int_{z_0}^{z_0+\omega_1}\frac{f'(z)}{f(z)}dz,\nonumber\\
&=&-\omega_2\ln\frac{f(z_0+\omega_1)}{f(z_0)},\nonumber\\
&=&2\pi in'\omega_2,\quad n'\in \mathbb{Z}.\nonumber
\end{eqnarray}
De même, on trouve
$$\int_{\gamma_2} z\frac{f'(z)}{f(z)}dz+\int_{\gamma_4}
\zeta\frac{f'(\zeta)}{f(\zeta)}d\zeta=2\pi in\omega_2,\quad n\in
\mathbb{Z}.$$ Par conséquent
\begin{eqnarray}
\sum_{k=1}^l n_ka_k-\sum_{k=1}^m
p_kb_k&=&n\omega_1+n'\omega_2,\nonumber\\
&=&\omega,\nonumber\\
&=&\mbox{période}.\nonumber
\end{eqnarray} $\square$

\begin{prop}
Soient $f$ et $g$ deux fonctions elliptiques ayant mêmes périodes.
Alors, il existe une relation algébrique de la forme
$$P(f(z) ,g(z))=0,$$
où $P$ est un polynôme à deux indéterminées et à coefficients
constants.
\end{prop}
\emph{Démonstration}: Soient $a_k$, $1\leq k\leq m$, les points du
parallélogramme fondamental en lesquels $f$ et (ou) $g$ ont des
pôles d'ordre maximum $p_k$, $1\leq k\leq m$. Soit $Q(Z,W)$ un
polynôme sans terme constant, de degré $n$ par rapport à $Z$ et
$W$. L'idée de la preuve est la suivante : On construit le
polynôme $Q$ de telle fàçon que les hypothèses de la proposition
1.4 concernant la fonction
$$F(z)=Q[f(z), g(z)],$$
soient satisfaites. La fonction $F(z)$ se réduit donc à une
constante $C$ et il suffit de choisir $P=Q-C$. En effet, la
fonction $F(z)$ est elliptique avec les mêmes périodes que les
fonctions $f(z)$, $g(z)$ et ne peut admettre de pôles qu'aux
points $a_k$. Les développements des fonctions $f$ et $g$ en
séries de Laurent au voisinage de $a_k$ ne contiennent que des
termes en $\frac{1}{(z-a_k)^j}$ avec $j\leq p_k$. La fonction
$F(z)$ ne peut avoir des pôles qu'aux points $a_k$ et son
développement en série de Laurent au voisinage de $a_k$ ne
contient que des termes en $\frac{1}{(z-a_k)^j}$ avec $j\leq p$ où
$n$ est le degré du polynôme $Q$ et $p=p_1+...+p_m$ est la somme
des ordres maximaux des fonctions $f, g$ aux points $a_k$. On
choisit les coefficients du polynôme $Q$ de manière à ce que les
parties principales de son développement en série de Laurent au
voisinage de $a_k$ soient nulles. Autrement dit, de sorte que le
développement en question ne contient pas des termes en
$\frac{1}{(z-a_k)^j}$ avec $j\leq p$. Donc l'élimination des pôles
de la fonction $F(z)$ fournira un système homogène de $np$
équations linéaires par rapport aux coefficients du polynôme $Q$.
Ce dernier étant de degré $n$ et comme il est supposé sans terme
constant, on aura donc $\frac{n(n+3)}{2}$ coefficients. En prenant
$\frac{n(n+3)}{2}>np$, on en déduit que le nombre des coefficients
(inconnues) est supérieur à celui des équations. Par conséquent,
le système en question admet au moins une solution non triviale
(i.e., non nulle). Finalement d'après la proposition 1.4, la
fonction $F(z)=Q[f(z), g(z)]$, est une constante $C$ et il suffit
de choisir $P=Q-C$. $\square$

\begin{cor}
Toute fonction elliptique $f(z)$ satisfait à une équation
différentielle de la forme
$$P(f(z) ,f'(z))=0,$$
où $P$ est un polynôme à deux indéterminées et à coefficients
constants.
\end{cor}
\emph{Démonstration}: D'après la remarque 1.5, la dérivée $f'(z)$
de la fonction elliptique $f(z)$ est aussi une fonction elliptique
et il suffit de poser $g(z)=f'(z)$ dans la proposition précédente.
$\square$

\section{Fonctions de Weierstrass}

Dans cette section on étudiera tout d'abord la fonction $\wp$ de
Weierstrass; c'est une fonction elliptique d'ordre 2 qui a un pôle
double à l'origine en tout point du parallélogramme fondamental.
Ensuite on introduit les deux autres fonctions de Weierstrass : la
fonction $\zeta$ et la fonction $\sigma$. Contrairement à la
fonction $\wp$, la fonction $\zeta$ est une fonction méromorphe
avec un pôle simple dans le parallélogramme fondamental tandis que
la fonction $\zeta$ est une fonction holomorphe partout. Les
fonctions de Weierstrass interviennent souvent lors de la
résolution de problèmes théoriques.

\subsection{Fonction $\wp$ de Weierstrass}

La fonction $\wp $ de Weierstrass est définie par
\begin{equation}\label{eqn:euler}
\wp(z)=\frac{1}{z^2}+\sum_{\omega \in \Lambda \backslash
\{0\}}(\frac{1}{(z-\omega )^2} -\frac{1}{\omega ^2}),
\end{equation}
où $\Lambda =\mathbb{Z}\omega _{1}\oplus\mathbb{Z}\omega _{2},$
est le réseau  engendré par deux nombres complexes $\omega _{1}$\
et $\omega _{2}$ diff\'{e}rents de $0$ tels que: $\mbox{Im}\left(
\frac{\omega _{2}}{\omega _{1}}\right)>0$.

\begin{prop}
La série (2.1) converge normalement sur tout compact ne
rencontrant pas le réseau $\Lambda$, i.e., sur tout compact de
$\mathbb{C}\backslash \Lambda$.
\end{prop}
\emph{Démonstration}: On montre que la série converge normalement
sur tout disque compact $\{z:\quad \mid z\mid\leq r\}$. Notons que
tout disque fermé ne contient qu'un nombre fini
d'éléments\footnote{Notons que $\mid x\omega_1+y\omega_2\mid$ est
une norme sur $\mathbb{R}^2$. Comme elle est équivalente à
$\sqrt{x^2+y^2}$, on peut donc trouver un $c>0$ tel que : $\mid
m\omega_1+n\omega_2\mid \geq c\sqrt{m^2+n^2}$, $\forall m,n$.} de
$\Lambda$ et que la nature de la série ne change évidemment pas si
on enlève ces éléments. Pour $\mid \omega \mid$ suffisamment
grand, on choisit $\mid \omega \mid\geq 2r$ pour tous les $\omega$
sauf un nombre fini; ceux qui sont dans le disque. On a
$$\mid \frac{1}{(z-\omega)^2}-\frac{1}{\omega^2}\mid=
\frac{\mid z \mid}{\mid \omega^3 \mid}.\frac{\mid
2-\frac{z}{\omega} \mid}{\mid 1-\frac{z}{\omega} \mid^2}.$$
Or
\begin{eqnarray}
\mid z \mid&\leq& r,\nonumber\\
\mid 2-\frac{z}{\omega} \mid&\leq& \frac{5}{2},\nonumber\\
\mid 1-\frac{z}{\omega} \mid&\geq& \frac{1}{2},\nonumber
\end{eqnarray}
donc
$$\mid \frac{1}{(z-\omega)^2}-\frac{1}{\omega^2}\mid\leq \frac{10 r}{\mid \omega
\mid^3},$$ et il suffit de prouver que la série
$$\sum_{\omega\in
\Lambda\backslash \{0\}}\frac{1}{\mid \omega^3 \mid},$$ converge.
Pour celà, considérons le parallélogramme
$$\Lambda_k=\{x\omega_1+y\omega_2 : \sup\{\mid x\mid, \mid
y\mid\}=k\},$$ où $n$ est un entier. Sur le parallélogramme
$\Lambda_1$  de cotés $2\omega_1$ et $2\omega_2$ dont le centre
est $0$, il y a $8$ points de $\Lambda$. Soit $d$ la plus courte
distance du point $z=0$ aux points de $\Lambda_1$. Pour chacun de
ces $8$ points, la distance à $0$ est $\geq d$, d'où
$$\frac{1}{\mid \omega\mid^3}\leq\frac{1}{d^3},$$ et
$$\sum_{\omega\in
\Lambda_1\backslash \{0\}}\frac{1}{\mid
\omega\mid^3}\leq\frac{8}{d^3}.$$ Sur le parallélogramme
$\Lambda_2$ (image de $\Lambda_1$ dans l'homothétie de centre $0$,
de rapport $2$), il y a $8\times 2=16$ points de $\Lambda$. Soit
$2d$ la plus courte distance du point $z=0$ aux points de
$\Lambda_2$. Pour chacun de ces $8$ points, la distance à $0$ est
$\geq 2d$, d'où
$$\sum_{\omega\in
\Lambda_2\backslash \{0\}}\frac{1}{\mid
\omega\mid^3}\leq\frac{8\times 2}{2^3d^3}=\frac{8}{2^2d^3}.$$ En
général sur le parallélogramme $\Lambda_k$ (image de $\Lambda_1$
dans l'homothétie de centre $0$, de rapport $k$), il y a $8 k$
points de $\Lambda$ et pour chacun de ces points, la distance à
$0$ est $\geq kd$. Dès lors,
$$\sum_{\omega\in
\Lambda_k\backslash \{0\}}\frac{1}{\mid
\omega\mid^3}\leq\frac{8k}{k^3d^3}=\frac{8}{k^2d^3}.$$ Ainsi la
série $$\sum_{\omega\in \Lambda\backslash \{0\}}\frac{1}{\mid
\omega^3 \mid},$$ est majorée par la série convergente
$$\frac{8}{d^3}\sum_{k=1}^\infty \frac{1}{k^2},$$
et par conséquent elle converge aussi en vertu du critère de
comparaison. $\square$

\begin{prop}
$\wp(z)$ est une fonction elliptique de périodes $\omega _{1}$\ et
$\omega _{2}$. Elle est paire et admet des pôles doubles aux
points $\omega\in \Lambda$, dont le résidu est nul. En outre,
$\wp'(z)$ est une fonction doublement périodique et elle est
impaire.
\end{prop}
\emph{Démonstration}: Notons tout d'abord que la fonction $\wp(z)$
est paire :
\begin{eqnarray}
\wp(-z)&=&\frac{1}{z^2}+\sum_{\omega \in \Lambda \backslash
\{0\}}(\frac{1}{(z-(-\omega))^2} -\frac{1}{(-\omega)
^2}),\nonumber\\
&=&\wp(z),\nonumber
\end{eqnarray}
car il suffit de remplacer $\omega$ par
$-\omega$. La dérivée de la fonction $\wp(z)$ est
\begin{eqnarray}
\wp'(z)&=&-\frac{2}{z^3}-2\sum_{\omega \in \Lambda \backslash
\{0\}}(\frac{1}{(z-\omega)^3},\nonumber\\
&=&-2\sum_{\omega\in \Lambda}(\frac{1}{(z-\omega)^3}.\nonumber
\end{eqnarray}
La fonction $\wp'(z)$ est doublement périodique de périodes
$\omega_1$ et $\omega_2$. En effet, on a
\begin{eqnarray}
\wp'(z+\omega_1)&=&-2\sum_{\omega\in
\Lambda}(\frac{1}{(z-(\omega-\omega_1))^3},\nonumber\\
&=&\wp'(z),
\end{eqnarray}
car $\omega-\omega_1$ est aussi une période. De façon analogue, on
montre que
\begin{equation}\label{eqn:euler}
\wp'(z+\omega_2)=\wp'(z),
\end{equation}
et donc
$$\wp'(z+\omega)=\wp'(z),\quad \forall \omega \in \Lambda.$$
En outre la fonction $\wp'(z)$ est impaire :
$$\wp'(-z)=-\wp'(z).$$
Montrons maintenant que $\wp(z)$ est une fonction elliptique de
périodes $\omega _{1}$\ et $\omega _{2}$. En intégrant les
relations (2.2) et (2.3), on obtient
$$\wp(z+\omega_1)-\wp(z)=C_1,$$
et
$$\wp(z+\omega_2)-\wp(z)=C_2,$$
où $C_1$ et $C_2$ sont des constantes. Posons
$z=-\frac{\omega_1}{2}$ et $z=-\frac{\omega_2}{2}$ (rappelons que
ces points ne sont pas des pôles de $\wp(z)$) dans la première et
seconde équation respectivement :
$$\wp(\frac{\omega_1}{2})-\wp(-\frac{\omega_1}{2})=C_1,$$
$$\wp(\frac{\omega_2}{2})-\wp(-\frac{\omega_2}{2})=C_2.$$
Or on a vu ci-dessus que la fonction $\wp(z)$ est paire, donc
$C_1=C_2=0$ et par conséquent
$$\wp(z+\omega_1)=\wp(z),$$
$$\wp(z+\omega_2)=\wp(z),$$
i.e., la fonction $\wp(z)$ est doublement périodique de périodes
$\omega_1$ et $\omega_2$. D'après la proposition 2.1, la série
(2.1) de fonctions méromorphes converge normalement sur tout
compact de $\mathbb{C}\backslash \Lambda$ et par conséquent sa
somme $\wp(z)$ est une fonction méromorphe sur $\mathbb{C}$. On en
déduit que $\wp(z)$ est une fonction elliptique de périodes
$\omega _{1}$\ et $\omega _{2}$. Notons enfin qu'au voisinage de
$z=\omega$, on a
$$\wp(z)=\frac{1}{(z-\omega)^2}+\mbox{fonction holomorphe},$$
ce qui signifie que les points $\omega\in \Lambda$ sont des pôles
doubles dont le résidu est nul. $\square$

\begin{prop}
Le développement de $\wp(z)$ en série de Laurent au voisinage du
point 0 est donné par
$$\wp(z)=\frac{1}{z^2}+\sum_{k=1}^\infty (2k+1)G_{2k+2}z^{2k},$$
où
$$G_k\equiv G_k(\Lambda)=\sum_{\omega \in \Lambda \backslash
\{0\}}\frac{1}{\omega ^k},\quad k\geq 4$$ et $G_k=0$ pour $k$
impaire.
\end{prop}
\emph{Démonstration}: Au voisinage de $z=0$, on a
$$\wp(z)=\frac{1}{z^2}+f(z),$$
où $f(z)$ est une fonction holomorphe. On a
$$f(z)=\sum_{\omega \in \Lambda \backslash
\{0\}}(\frac{1}{(z-\omega )^2} -\frac{1}{\omega ^2}),$$ avec
$f(0)=0$. Comme $\wp(z)$ est une fonction paire, alors au
voisinage de $z=0$ le développement de $f(z)$ en série de Laurent
a la forme
$$f(z)=a_2z^2+a_4z^4+...+a_{2k}z^{2k}+...,$$
avec
\begin{eqnarray}
a_2&=&\frac{g''(0)}{2}=3\sum_{\omega \in \Lambda \backslash
\{0\}}\frac{1}{\omega^4},\nonumber\\
a_4&=&\frac{g^{(4)}(0)}{4!}=5\sum_{\omega \in \Lambda \backslash
\{0\}}\frac{1}{\omega^6},\nonumber\\
&\vdots&\nonumber\\
a_{2k}&=&\frac{g^{(2k)}(0)}{(2k)!}=(2k+1)!\sum_{\omega \in \Lambda
\backslash \{0\}}\frac{1}{\omega^{2k+2}},\nonumber
\end{eqnarray}
Donc au voisinage de $z=0$, $\wp(z)$ admet un développement en
série de Laurent :
$$\wp(z)=\frac{1}{z^2}+\sum_{k=1}^\infty (2k+1)G_{2k+2}z^{2k},$$
avec
$$G_k\equiv G_k(\Lambda)=\sum_{\omega \in \Lambda \backslash
\{0\}}\frac{1}{\omega ^k},\quad k\geq 4.$$ $\square$

\begin{rem}
Nous allons donner une autre preuve similaire à la précédente.
Rappelons que :
$$\frac{1}{1-z}=\sum_{k=0}^\infty z^k,\quad \mid z \mid<1,$$
et
\begin{eqnarray}
\frac{1}{(1-z)^2}&=&(\frac{1}{1-z})',\nonumber\\
&=&\sum_{k=1}^\infty kz^{k-1},\nonumber\\
&=&\sum_{k=0}^\infty (k+1)z^k,\quad \mid z \mid<1.,\nonumber
\end{eqnarray}
Donc pour $\mid z \mid<\omega$, on a
\begin{eqnarray}
\frac{1}{(z-\omega)^2}-\frac{1}{\omega^2}&=&
\frac{1}{\omega^2}[\frac{1}{(1-\frac{z}{\omega})^2}-1],\nonumber\\
&=&\sum_{k=1}^\infty \frac{k+1}{\omega^{k+2}}z^k.\nonumber
\end{eqnarray}
Dès lors
$$\wp(z)=\frac{1}{z^2}+\sum_{\omega \in \Lambda \backslash
\{0\}}\sum_{k=1}^\infty \frac{k+1}{\omega^{k+2}}z^k.$$ En tenant
compte du fait que la fonction $\wp(z)$ est paire et que cette
double série est absolument convergente, on obtient
$$\wp(z)=\frac{1}{z^2}+3z^2\sum_{\omega \in \Lambda \backslash
\{0\}}\frac{1}{\omega^4}+5z^4\sum_{\omega \in \Lambda \backslash
\{0\}}\frac{1}{\omega^6}+...$$
\end{rem}

\begin{defn}
Les deux nombres complexes $g_2$ et $g_3$ définis par les séries
(dites d'Eisenstein) :
\begin{eqnarray}
g_{2}&=&60\sum_{\omega \in \Lambda \backslash \{0\}}
\frac{1}{\omega ^{4}},\nonumber\\
g_{3}&=&140\sum_{\omega \in \Lambda \backslash
\{0\}}\frac{1}{\omega^{6}},\nonumber
\end{eqnarray}
s'appellent invariants (de Weierstrass) de la fonction $\wp(z)$.
\end{defn}

\begin{prop}
La fonction $\wp(z)$ est solution dans $\Lambda $ de l'équation
différentielle:
\begin{equation}\label{eqn:euler}
(\wp ^{\prime }(z))^2=4(\wp(z))^3-g_{2}\wp(z)-g_3,
\end{equation}
où $g_{2}$ et $g_{3}$ sont les invariants de la fonction $\wp(z)$.
\end{prop}
\emph{Démonstration}: En utilisant les notations $g_2$ et $g_3$
introduites dans la définition 2.4, on récrit la fonction $\wp(z)$
sous la forme
$$\wp(z)=\frac{1}{z^2}+\frac{g_2}{20}z^2+\frac{g_3}{28}z^4+...$$
Les conditions de dérivation terme à terme de cette série étant
satisfaites, on obtient
$$\wp'(z)=\frac{-2}{z^3}+\frac{g_2}{10}z+\frac{g_3}{7}z^3+...$$
En élevant $\wp'(z)$ au carré et $\wp(z)$ au cube, on obtient
$$(\wp'(z))^2=\frac{4}{z^6}(1-\frac{g_2}{10}z^4-\frac{g_3}{7}z^6+...),$$
et
$$(\wp(z))^3=\frac{1}{z^6}(1+\frac{3g_2}{20}z^4-\frac{3g_3}{28}z^6+...).$$
Dès lors,
$$(\wp'(z))^2-4\wp(z))^3+g_2\wp(z)=-g_3+\frac{g_2^2}{20}z^2+\frac{g_2g_3}{28}z^4+...$$
La fonction
$$(\wp'(z))^2-4\wp(z))^3+g_2\wp(z)+g_3=\frac{g_2^2}{20}z^2+\frac{g_2g_3}{28}z^4+...$$
est holomorphe au voisinage de $z=0$ et elle est nulle en ce
point. Or cette fonction est doublement périodique, donc elle est
holomorphe au voisinage de tout point du parallélogramme
fondamentale et par conséquent elle est holomorphe dans tout
$\mathbb{C}$. Comme elle n'a pas de pôles, elle est bornée dans le
parallélogramme fondamentale (un compact) et donc bornée dans
$\mathbb{C}$. D'après le théorème de Liouville, cette fonction est
constante et puisqu'elle est nulle en $0$, elle est donc
identiquement nulle. $\square$

\begin{prop}
La fonction $\wp'(z)$ sannule au point $a\in \mathbb{C}$ tel que :
$-a\equiv a$ mod. $\Lambda$, i.e, $2a\in \Lambda$, $a\notin
\Lambda$. Autrement dit, modulo $\Lambda$, $\wp'(z)$ a trois zéros
simples : $\frac{\omega_1}{2}$, $\frac{\omega_2}{2}$,
$\frac{\omega_1+\omega_2}{2}$. En outre, en posant
$$e_1=\wp(\frac{\omega_1}{2}),\quad e_2=\wp(\frac{\omega_2}{2}),\quad
e_3=\wp(\frac{\omega_1+\omega_2}{2}),$$ on obtient
$$e_1\neq e_2\neq e_3,$$
et
$$
\left\{\begin{array}{rl}
e_1+e_2+e_3&=0,\\
e_1e_2+e_2e_3+e_3e_1&=-\frac{g_2}{4},\\
e_1e_2e_3&=\frac{g_3}{4}.
\end{array}\right.
$$
\end{prop}
\emph{Démonstration}: En tenant compte du fait que la fonction
$\wp'(z)$ est impaire et qu'elle est doublement périodique, on
obtient
\begin{eqnarray}
\wp'(\frac{\omega_k}{2})&=&-\wp'(-\frac{\omega_k}{2}),\nonumber\\
&=&-\wp'(-\frac{\omega_k}{2}+\omega_k),\nonumber\\
&=&-\wp'(\frac{\omega_k}{2}),\nonumber
\end{eqnarray}
où $k=1 ,2$ et donc $\wp'(\frac{\omega_k}{2})=0$. De même, on a
\begin{eqnarray}
\wp'(\frac{\omega_1+\omega_2}{2})&=&-\wp'(-\frac{\omega_1+\omega_2}{2}),\nonumber\\
&=&-\wp'(-\frac{\omega_1+\omega_2}{2}+\omega_1+\omega_2),\nonumber\\
&=&-\wp'(\frac{\omega_1+\omega_2}{2}),\nonumber
\end{eqnarray}
et donc $\wp'(\frac{\omega_1+\omega_2}{2})=0$. D'après la
proposition 2.5, on a
\begin{eqnarray}
(\wp ^{\prime }(z))^2&=&4(\wp(z))^3-g_{2}\wp(z)-g_3,\nonumber\\
&=&4(\wp(z)-e_1)(\wp(z)-e_2)(\wp(z)-e_3).\nonumber
\end{eqnarray}
Montrons que : $e_1\neq e_2\neq e_3$. En effet, on a vu que la
fonction $\wp(z)-e_j$, $(j=1 ,2, 3)$, est elliptique, possède un
pôle double et un zéro double. Dès lors, les relations :
$\wp(\frac{\omega_1}{2})-e_1=0$ et $\wp'(\frac{\omega_1}{2})=0$
signifient que la fonction $\wp(z)-e_1$ a un zéro double en
$\frac{\omega_1}{2}$ ce qui implique que : $e_1\neq e_2$ et
$e_1\neq e_3$ car sinon la fonction $\wp(z)$ aurait plus d'un zéro
double, ce qui contredit la multiplicité. Pour les points
$\frac{\omega_2}{2}$ et $\frac{\omega_1+\omega_2}{2}$ il suffit de
faire un raisonnement similaire au précédent. Par conséquent,
$e_1\neq e_2\neq e_3$. Les autres relations entre les coefficients
de l'équation (2.4) et ses racines, découlent immédiatemment des
propriétés des racines des équations algébriques. $\square$

\begin{rem}
En posant $w=\wp(z)$, l'équation (2.4) s'écrit
$$(\frac{dw}{dz})^2=4w^3-g_{2}w-g_3.$$
Or $z\rightarrow 0$ lorsque $w\rightarrow \infty$, donc
\begin{equation}\label{eqn:euler}
z=\int_\infty^w\frac{dw}{\sqrt{4w^3-g_{2}w-g_3}}.
\end{equation}
Autremant dit, la fonction $w=\wp(z)$ s'obtient par inversion de
l'intégrale (2.5) (dite intégrale elliptique sous forme de
Weierstrass). Réciproquement, si le polynôme $4w^3-g_{2}w-g_3$ n'a
pas de zéros multiples (i.e., son discriminant est non nul :
$g_2^3-27g_3^2\neq 0$), alors l'inversion de l'intégale (2.5)
conduit à la fonction $\wp(z)$ de Weierstrass.
\end{rem}

\begin{prop}
L'application
\begin{eqnarray}
\mathbb{C}/\Lambda \longrightarrow \mathbb{CP}^{2},&&z\longmapsto
[1,\wp(z),\wp '(z)], z\neq 0,\nonumber\\
&&0\longmapsto [0,0,1],\nonumber
\end{eqnarray}
est un isomorphisme entre le tore complexe $\mathbb{C}/\Lambda $
et la courbe elliptique $\mathcal{E}$ d'\'{e}quation affine:
\begin{equation}\label{eqn:euler}
y^{2}=4x^{3}-g_{2}x-g_{3}.
\end{equation}
\end{prop}
\emph{Démonstration}: Il suffit de poser $x=\wp(z)$, $y=\wp'(z)$
et d'utiliser l'équation différentielle (2.4). $\square$

\begin{prop}
Soient $u, v\notin\Lambda$ et $u\pm v\notin\Lambda$, alors la
fonction $\wp(z)$ vérifie la loi d'addition
$$\wp(u)+\wp(v)+\wp(u+v)
=\frac{1}{4}(\frac{\wp^{\prime}(u)-\wp^{\prime}(v) }{\wp(u)-\wp (
v)})^2,$$ ainsi que la formule de duplication
$$\wp(2z)=\frac{1}{4}(\frac{\wp^{\prime \prime }
(z)}{\wp^{\prime}(z)})^{2}-2\wp(z).$$ Géométriquement, si $P_1$ et
$P_2$ sont deux points distincts de la courbe elliptique
$\mathcal{E}$ d'équation affine (2.6), alors $\mathcal{E}$ a trois
points d'intersection avec la droite sécante $L=P_1P_2$ qui passe
par $P_1$ et $P_2$. Dans le cas où $P_1=P_2$, alors $L$ est la
tangente à la courbe au point $P_1=P_2$. Les coordonnées du
$3^{\mbox{ème}}$ point d'intersection s'expriment comme fonctions
rationnelles des deux autres.
\end{prop}
\emph{Démonstration}: Considérons la courbe elliptique
$\mathcal{E}$ d'\'{e}quation affine
$$y^{2}=4x^{3}-g_{2}x-g_{3}.$$
Soient $P_1$ et $P_2$ sont deux points distincts de $\mathcal{E}$,
ayant pour coordonnées $(x_1=\wp(u), y_1=\wp'(u))$ et
$(x_2=\wp(v), y_2=\wp'(v))$ respectivement. Soit $y=ax+b$
l'équation de la droite sécante $L=P_1P_2$ qui passe par $P_1$ et
$P_2$. Cette droite coupe la courbe $\mathcal{E}$ de telle façon
que :
$$y^2=4x^3-g_2x-g_3=(ax+b)^2,$$
ou ce qui revient au même
$$\varphi(x)\equiv 4x^3-g_2x-g_3-(ax+b)^2=0.$$
Considérons maintenant la fonction elliptique
$$f(z)=\wp'(z)-a\wp(z)-b.$$
Comme $\wp'(z)$ a un pôle d'ordre 3 à l'origine, il en est donc de
même pour $f(z)$. Cette fonction a donc trois zéros dont deux sont
connus : $z=u$, $z=v$ et un troisième que nous noterons
provisoirement $z=t$. D'après la proposition 1.7, la somme des
pôles est égale à celui des zéros, d'où
$$0+0+0=u+v+t,\quad (\mbox{mod.} \Lambda),$$
et donc $t=-u-v$. Dès lors,
\begin{eqnarray}
0&=&f(u)=\wp'(u)-a\wp(u)-b,\nonumber\\
0&=&f(v)=\wp'(v)-a\wp(v)-b,\\
0&=&f(-u-v)=-\wp'(u+v)-a\wp(u+v)-b,\nonumber
\end{eqnarray}
en tenant compte du fait que $\wp(z)$ est paire et $\wp'(z)$ est
impaire. On déduit immédiatement des deux premières équations que
:
\begin{eqnarray}
a&=&\frac{\wp'(u)-\wp'(v)}{\wp(u)-\wp(v)},\nonumber\\
b&=&\frac{\wp'(v)\wp(u)-\wp'(u)\wp(v)}{\wp(u)-\wp(v}.\nonumber
\end{eqnarray}
En remplaçant ces expressions dans la troisième équation du
système (2.7), on obtient la relation
\begin{equation}\label{eqn:euler}
\wp'(u+v)=-\frac{\wp'(u)-\wp'(v)}{\wp(u)-\wp(v)}\wp(u+v)
-\frac{\wp'(v)\wp(u)-\wp'(u)\wp(v)}{\wp(u)-\wp(v)}.
\end{equation}
Par ailleurs, on a
\begin{eqnarray}
\varphi(\wp(z))&=&4\wp^3(z)-g_2\wp(z)-g_3-(a\wp(z)+b)^2,\nonumber\\
&=&4\wp^3(z)--a^2\wp^2(z)-(g_2+2ab)\wp(z)-g_3-b^2,\nonumber
\end{eqnarray}
et puisque
$$\varphi(\wp(u))=\varphi(\wp(v))=\varphi(\wp(u+v))=0,$$
alors
$$\wp(u)+\wp(v)+\wp(u+v)=\frac{a^2}{4}.$$
En remplaçant $a$ par sa valeur obtenue précédemment, on obtient
\begin{equation}\label{eqn:euler}
\wp(u)+\wp(v)+\wp(u+v)=\frac{1}{4}(\frac{\wp^{\prime}(u)-\wp^{\prime}(v)}{\wp(u)-\wp(v)})^2.
\end{equation}
Rappelons que la courbe $\mathcal{E}$ a deux points d'intersection
$P_1$ de coordonnées $(x_1=\wp(u), y_1=\wp'(u))$ et $P_2$ de
coordonnées $(x_2=\wp(v), y_2=\wp'(v))$ avec la droite sécante
$L=P_1P_2$ passant par $P_1$ et $P_2$. On sait qu'il existe un
troisième point unique $P_3\in\mathcal{E}\cap L$ de coordonnées
$(x_3=\wp(w), y_3=\wp'(w))$. D'après le système (2.7), les
coordonnées $(x_3 ,y_3)$ s'expriment en fonction de $(x_1 ,y_1)$
et $(x_2 ,y_2)$ comme suit
\begin{eqnarray}
x_3&=&-(x_1+x_2)+\frac{1}{4}(\frac{y_1-y_2}{x_1-x_2})^2,\nonumber\\
y_3&=&ax_3+b,\nonumber\\
&=&(\frac{y_1-y_2}{x_1-x_2})[-(x_1+x_2)+\frac{1}{4}(\frac{y_1-y_2}{x_1-x_2})^2]
+\frac{y_2x_1-y_1x_2}{x_1-x_2}.\nonumber
\end{eqnarray}
Dans la formule (2.9), divisons le numérateur et le dénominateur
par $u-v$,
$$\wp(u)+\wp(v)+\wp(u+v)
=\frac{1}{4}[\frac{\frac{\wp'(u)-\wp'(v)}{u-v}}{\frac{\wp(u)-\wp(v)}{u-v}}
]^2.$$ En faisant tendre $u$ et $v$ vers $z$, on obtient
$$\wp(2z)=\frac{1}{4}(\frac{\wp^{\prime \prime }
(z)}{\wp^{\prime}(z)})^{2}-2\wp(z),$$ et les coordonnées $(x_3,
y_3)$ du troisième point $P_3\in\mathcal{E}\cap L$ deviennent
\begin{eqnarray}
x_3&=&-2x_1+\frac{1}{4}(\frac{12x_1^2-g_2}{2y_1})^2,\nonumber\\
y_3&=&-y_1+\frac{1}{4}(\frac{12x_1^2-g_2}{2y_1})(x_1-x_3).\nonumber
\end{eqnarray} $\square$

\begin{rem}
Notons que nous avons choisi $a$ et $b$ de telle façon que :
\begin{eqnarray}
0&=&f(u)=\wp'(u)-a\wp(u)-b,\nonumber\\
0&=&f(v)=\wp'(v)-a\wp(v)-b.\nonumber
\end{eqnarray}
Il faut donc que
$$
\det\left(\begin{array}{cc}
\wp(u)&1\\
\wp(v)&1
\end{array}\right)
=\wp(u)-\wp(v)\neq0.
$$
Evidemment si $\wp(u)-\wp(v)=0$, alors il suffit de déplacer $u$
et $v$ légérement de façon à avoir $\wp(u)-\wp(v)\neq0$. Le
système (2.7) s'écrit
$$
\left(\begin{array}{ccc}
\wp'(u)&\wp(u)&1\\
\wp'(v)&\wp(v)&1\\
-\wp'(u+v)&\wp(u+v)&1
\end{array}\right)
\left(\begin{array}{c}
1\\
a\\
b
\end{array}\right)
=0,
$$
et comme le déterminant ci-dessus est $\neq0$, alors on obtient la
condition
$$
\det\left(\begin{array}{ccc}
\wp'(u)&\wp(u)&1\\
\wp'(v)&\wp(v)&1\\
-\wp'(u+v)&\wp(u+v)&1
\end{array}\right)=0,
$$
i.e., la relation (2.8) obtenue précédemment.
\end{rem}

Soit $\mathcal{E}_\Lambda$ l'ensemble des fonctions elliptiques.
Cet ensemble est un espace vectoriel (et même un corps). On note
$\mathbb{C}(X)$ l'ensemble des fonctions rationnelles d'une
variable.

\begin{prop}
On a $\mathcal{E}_\Lambda=\mathbb{C}(\wp, \wp'),$ i.e., toute
fonction elliptique pour $\Lambda$ est une fonction rationnelle de
$\wp(z)$ et $\wp'(z)$. Plus précisement, l'application
$$\mathbb{C}(X)\times\mathbb{C}(X)\longrightarrow \mathcal{E}_\Lambda,\quad
(g,h)\longmapsto f(z)=g(\wp(z))+\wp'(z)h(\wp(z)),$$ est un
isomorphisme entre espaces vectoriels.
\end{prop}
\emph{Démonstration}: Soit $f\in \mathcal{E}_\Lambda$. On peut
évidemment écrire $f$ comme une somme d'une fonction paire et
d'une fonction impaire :
$$f(z)=\frac{f(z)+f(-z)}{2}+\frac{f(z)-f(-z)}{2}.$$
La fonction $\wp'(z)$ étant impaire, on réecrit la fonction $f$
sous la forme :
$$f(z)=\frac{f(z)+f(-z)}{2}+\wp'(z)(\frac{f(z)-f(-z)}{2\wp'(z)}).$$
Comme les fonctions $\frac{f(z)+f(-z)}{2}$ et
$\wp'(z)(\frac{f(z)-f(-z)}{2\wp'(z)})$ sont paires, il suffit donc
de démontrer que le sous-corps des fonctions elliptiques paires
par rapport à $\Lambda$ est engendré par $\wp(z)$. Soit donc $f$
une fonction elliptique paire telle que : $f\neq0$, $f\neq\infty$
aux points du parallélogramme des périodes (i.e., $f$ n'a ni pôle,
ni zéro sur le réseau). Si $z$ est un point tel que : $f(z)=0$,
alors comme $f$ est paire, $f(-z)=0$ et on aura un ordre pair. Dès
lors, on peut toujours choisir des points : $z_1,...,z_k,
-z_1,...,-z_k$ qui sont des zéros de $f$ et des points
$p_1,...,p_k, -p_1,...,-p_k$ qui sont des pôles de $f$.
Considérons la fonction
$$g(z)=\prod_{j=1}^k\frac{\wp(z)-\wp(z_j)}{\wp(z)-\wp(p_j)}.$$
La fonction $\wp(z)$ étant paire, alors les zéros (resp. pôles) de
$g(z)$ sont $z=z_j$ (resp. $p_j$) et $z=-z_j$ (resp. $-p_j$). La
fonction elliptique $g(z)$ a les mêmes pôles et les mêmes zéros
que $f(z)$. Dès lors, la fonction $\frac{g(z)}{f(z)}$ n'a pas de
pôles et n'a pas de zéros et d'après le théorème de Liouville elle
est constante. Par conséquent,
\begin{eqnarray}
f(z)&=&Cg(z),\quad (C=\mbox{constante}),\nonumber\\
&=&C\prod_{j=1}^k\frac{\wp(z)-\wp(z_j)}{\wp(z)-\wp(p_j)},\nonumber\\
&=&\mbox{fonction rationnelle de } \wp(z).\nonumber
\end{eqnarray}
Notons que si $f$ a un pôle ou un zéro dans le parallélogramme des
périodes, alors pour se débarasser du pôle ou du zéro, il suffit
de multiplier $f(z)$ par $(\wp(z))^j$. Autrement dit, la fonction
$f(z)(\wp(z))^j$ est paire, sans pôles, ni zéros en $(0, 0)$ et
c'est une fonction rationnelle de $\wp(z)$. $\square$

\begin{rem}
Posons $\mathcal{E}^+_\Lambda=\{f\in \mathcal{E}_\Lambda : f
\mbox{paire}\}$. Dans la preuve précédente, on a montré que :
$\mathcal{E}^+_\Lambda=\mathbb{C}(\wp),$ i.e., le sous-corps des
fonctions elliptiques paires par rapport à $\Lambda$ est engendré
par $\wp(z)$.
\end{rem}

\begin{rem}
D'après la proposition 2.9, pour caractériser le corps des
fonctions elliptiques on forme l'anneau quotient $\mathbb{C}[X,Y]$
par l'idéal principal correspondant à l'équation
$$Y^2=4X^3-g_2X-g_3.$$
Plus précisement, on a
$$\mathcal{E}_\Lambda=\mathbb{C}(\wp, \wp')\simeq
\mathbb{C}[X,Y]/(Y^2-4X^3+g_2X+g_3),$$ où $\wp=X$ et $\wp'$ est
identifiée à l'image de $Y$ dans le quotient.
\end{rem}

\subsection{Fonction $\zeta$\ de Weierstrass}

La fonction $\zeta$\ de Weierstrass\footnote{A ne pas confondre
avec la fonction $\zeta$ de Riemann.} est définie par
\begin{equation}\label{eqn:euler}
\zeta(z)=\frac{1}{z}-\int_0^z(\wp(z)-\frac{1}{z^2})dz.
\end{equation}
Notons que la dérivée de cette fonction est
\begin{equation}\label{eqn:euler}
\zeta'(z)=-\wp(z).
\end{equation}
En remplaçant $\wp(z)$ par (2.1) et après intégration on obtient
\begin{equation}\label{eqn:euler}
\zeta(z)=\frac{1}{z}+\sum_{\omega\in\Lambda\backslash\{0\}}
(\frac{1}{z-\omega}+\frac{1}{\omega}+\frac{z}{\omega^2}).
\end{equation}

\begin{prop}
a) La fonction $\zeta(z)$ est impaire.\\
b) $\zeta(z)$ n'est pas une fonction elliptique.\\
c) La fonction $\zeta(z)$ n'est pas périodique et on a
\begin{equation}\label{eqn:euler}
\zeta(z+\omega_k)-\zeta(z)=\tau_k,\quad(k=1 ,2)
\end{equation}
où $\tau_k$ sont des constantes.\\
d) Les nombres $\omega_k$ et $\tau_k$ sont liés par la relation de
Legendre : $$\tau_1\omega_2-\tau_2\omega_1=2\pi i.$$
\end{prop}
\emph{Démonstration}: a) En utilisant (2.11) et la parité de
$\wp(z)$, on obtient
\begin{eqnarray}
(\zeta(z)+\zeta(-z))'&=&\zeta'(z)-\zeta'(-z),\nonumber\\
&=&-\wp(z)+\wp(z),\nonumber\\
&=&0.\nonumber
\end{eqnarray}
D'où
$$\zeta(z)+\zeta(-z)=C,$$
où $C$ est une constante. En remplaçant $\zeta(z)$ par son
expression (2.10), on obtient
$$\int_{-z}^z(\wp(z)-\frac{1}{z^2})dz=C,$$
et
$$0=\lim_{z\rightarrow 0}\int_{-z}^z(\wp(z)-\frac{1}{z^2})dz=C.$$
Par conséquent $\zeta(z)=-\zeta(-z)$, i.e., la fonction $\zeta(z)$
est impaire.\\
b) En effet, $\zeta(z)$ a des pôles simples en $\omega$ et d'après
la remarque 1.4, il n'existe pas de fonction elliptique de premier
ordre.\\
c) D'après (2.11) et le fait que $\wp(z)$ est doublement
périodique, on a
\begin{eqnarray}
(\zeta(z+\omega_k)-\zeta(z))'&=&-\wp(z+\omega_k)+\wp(z),\nonumber\\
&=&-\wp(z)+\wp(z),\nonumber\\
&=&0,\nonumber
\end{eqnarray}
et par conséquent $\zeta(z+\omega_k)-\zeta(z)=\tau_k, (k=1 ,2)$
où $\tau_k$ sont des constantes.\\
d) Comme dans la preuve de la proposition 1.6, soit
$\gamma=\gamma_1\cup\gamma_2\cup\gamma_3\cup\gamma_4$ la frontière
du parallélogramme fondamental $\Pi$ relativement au réseau
$\Lambda$, avec $ \gamma_1=[z_0,z_0+\omega_1],
\gamma_2=[z_0+\omega_1,z_0+\omega_1+\omega_2],
\gamma_3=[z_0+\omega_1+\omega_2,z_0+\omega_2]$ et $
\gamma_4=[z_0+\omega_2,z_0]$. Supposons que l'unique pôle $z=0$ de
$\zeta(z)$ soit à l'intérieur de ce parallélogramme, sinon on peut
toujours en vertu de la proposition 1.5 choisir un autre
parallélogramme proche du précédent de façon à ce que le pôle en
question soit à son intérieur. Le résidu de $\zeta(z)$ au point
$z=0$ étant égal à 1, on déduit du théorème des résidus que :
\begin{equation}\label{eqn:euler}
\sum_{j=1}^4\int_{\gamma_j}\zeta(z)dz=2\pi i.
\end{equation}
Notons que
\begin{eqnarray}
\int_{\gamma_2}\zeta(z)dz&=&\int_{[z_0+\omega_1,z_0+\omega_1+\omega_2]}\zeta(z)dz,\nonumber\\
&=&\int_{[z_0,z_0+\omega_2]}\zeta(u+\omega_1)du,\quad u\equiv z-\omega_1,\nonumber\\
&=&\int_{[z_0,z_0+\omega_2]}\zeta(z+\omega_1)dz,\nonumber
\end{eqnarray}
et
\begin{eqnarray}
\int_{\gamma_3}\zeta(z)dz&=&\int_{[z_0+\omega_1+\omega_2,z_0+\omega_2]}\zeta(z)dz,\nonumber\\
&=&\int_{[z_0+\omega_1,z_0]}\zeta(v+\omega_2)du,\quad v\equiv z-\omega_2,\nonumber\\
&=&-\int_{[z_0,z_0+\omega_1]}\zeta(z+\omega_2)dz.\nonumber
\end{eqnarray}
En remplaçant ces expressions dans (2.14), on obtient
$$\int_{[z_0,z_0+\omega_1]}(\zeta(z)-\zeta(z+\omega_2))dz
+\int_{[z_0,z_0+\omega_2]}(\zeta(z+\omega_1)-\zeta(z))dz=2\pi i,$$
et d'après (2.13),
$$\int_{[z_0,z_0+\omega_1]}(-\tau_2)dz
+\int_{[z_0,z_0+\omega_2]}\tau_1dz=2\pi i,$$ i.e.,
$\tau_1\omega_2-\tau_2\omega_1=2\pi i$. $\square$

\subsection{Fonction $\sigma$\ de Weierstrass}

La fonction $\sigma$\ de Weierstrass est définie par
\begin{equation}\label{eqn:euler}
\sigma(z)=z e^{\int_0^z(\zeta(z)-\frac{1}{z})dz},
\end{equation}
et sa dérivée logarithmique est
\begin{equation}\label{eqn:euler}
(\ln \sigma(z))'=\frac{\sigma'(z)}{\sigma (z)}=\zeta(z).
\end{equation}
En remplaçant $\zeta(z)$ par son expression (2.12), on obtient
\begin{equation}\label{eqn:euler}
\sigma(z)=z \prod_{\omega\in
\Lambda\backslash\{0\}}(1-\frac{z}{\omega})e^{\frac{z}{\omega}
+\frac{1}{2}(\frac{z}{\omega})^2}.
\end{equation}

\begin{prop}
a) $\sigma(z)$ est une fonction impaire.\\
b) La fonction $\sigma(z)$ vérifie la relation
$$\sigma(z+\omega_k)=-e^{\tau_k(z+\frac{\omega_k}{2})}.\sigma(z),\quad
(k=1,2),$$ où $\tau_k$ sont des constantes.
\end{prop}
\emph{Démonstration}: a) D'après (2.17), on a
\begin{eqnarray}
\sigma(-z)&=&-z \prod_{\omega\in
\Lambda\backslash\{0\}}(1+\frac{z}{\omega})e^{-\frac{z}{\omega}
+\frac{1}{2}(\frac{z}{\omega})^2},\nonumber\\
&=&-z \prod_{\omega\in
\Lambda\backslash\{0\}}(1-\frac{z}{\eta})e^{\frac{z}{\eta}
+\frac{1}{2}(\frac{z}{\eta})^2},\quad \eta\equiv -\omega,\nonumber\\
&=&-\sigma(z).\nonumber
\end{eqnarray}
Une autre preuve consiste à utiliser l'autre définition (2.15) de
$\sigma(z)$. On a
\begin{eqnarray}
\sigma(-z)&=&-z e^{\int_0^{-z}(\zeta(u)-\frac{1}{u})du},\nonumber\\
&=&-z e^{-\int_0^{z}(\zeta(-v)+\frac{1}{v})dv},\quad v\equiv-u\nonumber\\
&=&-z e^{-\int_0^{z}(-\zeta(v)+\frac{1}{v})dv},\quad (\zeta \mbox{est impaire})\nonumber\\
&=&-z e^{\int_0^{z}(\zeta(v)-\frac{1}{v})dv},\nonumber\\
&=&-\sigma(z).\nonumber
\end{eqnarray}
b) On a
\begin{eqnarray}
\frac{\sigma'(z+\omega_k)}{\sigma(z+\omega_k)}&=&\zeta(z+\omega_k),
\quad\mbox{d'après} (2.16)\nonumber\\
&=&\zeta(z)+\tau_k,\quad\mbox{d'après} (2.13)\nonumber\\
&=&\frac{\sigma'(z)}{\sigma(z)}+\tau_k.\quad\mbox{d'après}(2.16)\nonumber
\end{eqnarray}
En intégrant, on obtient
$$\ln\sigma(z+\omega_k)=\ln\sigma(z)+\tau_kz
C_k,\quad C_k\equiv\mbox{constante}$$ d'où
$$\sigma(z+\omega_k)=e^{\tau_kz+C_k}.\sigma(z).$$
Pour $z=-\frac{\omega_k}{2}$, on a
$$\sigma(\frac{\omega_k}{2})=e^{-\frac{\tau_k\omega_k}{2}}e^{C_k}\sigma(-\frac{\omega_k}{2}).$$
Or $\sigma(z)$ est impaire, donc
$$e^{C_k}=-e^{\frac{\tau_k\omega_k}{2}},$$
et par conséquent
$\sigma(z+\omega_k)=-e^{\tau_k(z+\frac{\omega_k}{2})}.\sigma(z),\quad
(k=1,2)$. $\square$

\begin{prop}
Soit $f$ une fonction elliptique d'ordre $n$. Désignons par
$a_1,...,a_n$ (resp. $b_1,...,b_n$) les zéros (resp. pôles) de $f$
dans le parallélogramme des périodes. Ici tous les zéros et les
pôles sont comptés avec leurs ordres de multiplicités. Alors
$$f(z)=C\sigma(z+\sum_{j=2}^na_j-\sum_{j=1}^nb_j)\frac{\prod_{j=2}^n\sigma(z-a_j)}
{\prod_{j=1}^n\sigma(z-b_j)},$$ où $C$ est une constante.
\end{prop}
\emph{Démonstration}: D'après la proposition 1.8, on a
$$\sum_{j=1}^na_j-\sum_{j=1}^nb_j=\mbox{période}\equiv \omega.$$
D'où
$$\sum_{j=2}^na_j-\sum_{j=1}^nb_j=\omega-a_1.$$
Considérons la fonction
$$g(z)=\sigma(z+\omega-a_1)\frac{\prod_{j=2}^n\sigma(z-a_j)}
{\prod_{j=1}^n\sigma(z-b_j)}.$$ D'après la proposition 2.11 (point
b)), on a
$$\sigma(z+\omega_k)=-e^{\tau_k(z+\frac{\omega_k}{2})}.\sigma(z),\quad
(k=1,2),$$ et comme $\sigma(z)$ est une fonction impaire, alors
\begin{eqnarray}
g(z+\omega_k)&=&e^{\tau_k(\sum_{j=1}^nb_j-a_1+\omega-\sum_{j=2}^na_j)}
\sigma(z+a_1+\omega)\frac{\prod_{j=2}^n\sigma(z-a_j)}{\prod_{j=1}^n\sigma(z-b_j)},\nonumber\\
&=&g(z).\nonumber
\end{eqnarray}
La fonction $\frac{f(z)}{g(z)}$ n'a pas de pôles dans le
parallélogramme des périodes. Puisque cette fonction est
doublement périodique, alors elle est bornée sur $\mathbb{C}$ et
par conséquent, elle est constante en vertu du théorème de
Liouville. $\square$

\section{Intégrales elliptiques et fonctions de Jacobi}

Dans cette section, on va étudier les fonctions de Jacobi. Ce sont
des fonctions elliptiques du second ordre qui ont deux pôles
simples dans le parallélogramme des périodes. Ces fonctions
interviennent souvent lors de la résolution de problèmes
pratiques.

On appelle en général intégrale elliptique une intégrale de la
forme
$$\int R(s,\sqrt{P(s)})ds,$$
où $R$ est une fonction rationnelle à deux variables et $P(s)$ un
polynôme de degré 3 ou 4 avec des racines simples. En général,
cette intégrale ne s'exprime pas au moyen de fonctions
élémentaires c'est-à-dire celles que l'on obtient en appliquant à
la variable $s$ les opérations algébriques (addition, soustration,
multiplication, division) en nombre fini, ainsi que les fonctions
logarithmiques, trigonométriques et leurs inverses. Nous verrons
que les fonctions inverses de ces intégrales elliptiques sont des
fonctions elliptiques. On montre qu'à l'aide de transformations
élémentaires, une intégrale elliptique se ramène à l'une des
formes canoniques (de Legendre) :
\begin{eqnarray}
&&\int\frac{ds}{\sqrt{(1-s^2)(1-k^2s^2)}},\nonumber\\
&&\int\sqrt{\frac{1-k^2s^2}{1-s^2}}ds,\nonumber\\
&&\int\frac{ds}{(1+ls^2)\sqrt{(1-s^2)(1-k^2s^2)}},\nonumber
\end{eqnarray}
où $k$ et $l$ sont des constantes. La première de ces intégrales
est dite intégrale elliptique de première espèce, la seconde
intégrale elliptique de seconde espèce et la troisième intégrale
elliptique de troisième espèce. On peut écrire ces intégrales sous
une forme un peu différente, en posant $s=\sin \varphi,$ et les
intégrales précédentes s'écrivent
\begin{eqnarray}
&&\int\frac{d\varphi}{\sqrt{1-k^2\sin^2\varphi}},\nonumber\\
&&\int\sqrt{1-k^2\sin^2\varphi}d\varphi,\nonumber\\
&&\int\frac{d\varphi}{(1+l\sin^2\varphi)\sqrt{1-k^2\sin^2\varphi}}.\nonumber
\end{eqnarray}
En général, on adopte les notations suivantes :
\begin{eqnarray}
F(k,\varphi)&=&\int_0^{\sin\varphi}\frac{ds}{\sqrt{(1-s^2)(1-k^2s^2)}}
=\int_0^\varphi\frac{d\varphi}{\sqrt{1-k^2\sin^2\varphi}},\nonumber\\
E(k,\varphi)&=&\int_0^{\sin\varphi}\sqrt{\frac{1-k^2s^2}{1-s^2}}ds
=\int_0^\varphi\sqrt{1-k^2\sin^2\varphi}d\varphi,\nonumber\\
\Pi(k,l,\varphi)&=&\int_0^{\sin\varphi}\frac{ds}{(1+ls^2)\sqrt{(1-s^2)(1-k^2s^2)}}
=\int_0^\varphi\frac{d\varphi}{(1+l\sin^2\varphi)\sqrt{1-k^2\sin^2\varphi}}.\nonumber
\end{eqnarray}
On rencontre souvent des intégrales où la borne supérieure est
$\varphi=\frac{\Pi}{2}.$ Dans ce cas, on écrit
\begin{eqnarray}
F(k)&=&\int_0^1\frac{ds}{\sqrt{(1-s^2)(1-k^2s^2)}}
=\int_0^{\frac{\Pi}{2}}\frac{d\varphi}{\sqrt{1-k^2\sin^2\varphi}},\nonumber\\
E(k)&=&\int_0^1\sqrt{\frac{1-k^2s^2}{1-s^2}}ds
=\int_0^{\frac{\Pi}{2}}\sqrt{1-k^2\sin^2\varphi}d\varphi,\nonumber
\end{eqnarray}
et ces intégrales sont dites intégrales elliptiques complètes
respectivement de première et de seconde espèce. On montre que
\begin{eqnarray}
F(k)&=&\frac{\Pi}{2}(1+(\frac{1}{2})^2+(\frac{1.3}{2.4})^2k^4
+(\frac{1.3.5}{2.4.6})^2k^6+...),\nonumber\\
E(k)&=&\frac{\Pi}{2}(1-(\frac{1}{2})^2-(\frac{1.3}{2.4})^2\frac{k^4}{3}
-(\frac{1.3.5}{2.4.6})^2\frac{k^6}{5}-...).\nonumber
\end{eqnarray}
Considérons des intégrales elliptiques de la forme
$$
t=\int_0^s\frac{ds}{\sqrt{(1-s^2)(1-k^2s^2)}},\quad0\leq k\leq 1
$$
et voyons avec un peu plus de détail les propriétés de cette
intégrale de première espèce tout en sachant que les propriétés
des autres intégrales s'obtiennent de façon similaire. Nous avons
vu ci-dessus que le changement de variable $s=\sin \varphi,$
ramène cette intégrale à la forme
$$t=\int_0^\varphi\frac{d\varphi}{\sqrt{1-k^2\sin^2\varphi}} .$$
Nous envisagerons tout d'abord le cas où $k\neq 0$ et $k\neq1.$ La
fonction $t(\varphi)$ définie par cette intégrale est strictement
croissante et dérivable. Elle possède donc un inverse, qu'on
appelle amplitude de $t$ et qui se note
$$\varphi=\mathbf{am}t=\mathbf{am}(t;k).$$
Notons que si $k=0,$ alors
$$t=\int_0^s\frac{ds}{\sqrt{1-s^2}}=\arcsin s,$$
d'où $s=\sin t.$ Pour $k\neq 0,$ on note par analogie la fonction
inverse de l'intégrale en question par
$$s=\mathbf{sn}t=\mathbf{sn}(t;k),$$
que l'on nomme fonction elliptique de Jacobi (Lire s, n, t en
détachant les lettres). Le nombre $k$ est appelé module de la
fonction. Lorsqu'il n'y a pas ambiguité sur le module $k$,
on écrit tout simplement $\mathbf{sn}t$ au lieu de $\mathbf{sn}(t;k)$.\\
La fonction $\varphi=\mathbf{am} t$ est une fonction impaire
strictement croissante de $t$. Elle satisfait à
$$\mathbf{am} (0)=0,\qquad \frac{\partial\mathbf{am}}{\partial
t}(0)=1.$$ Comme $s=\sin \varphi$, on peut donc écrire
$s=\mathbf{sn} t=\sin (\mathbf{am} t).$\\
La deuxième et troisième fonction elliptique de Jacobi sont
définies respectivement par
$$\mathbf{cn} t=\mathbf{cn}(t;k)=\cos\mathbf{am}t,$$
et
$$\mathbf{dn} t=\mathbf{dn}(t;k)=\sqrt{1-k^2\mathbf{sn}^2 t}.$$
Pour $\mathbf{cn} t$, lire c, n, t en détachant les lettres. De
même, pour $\mathbf{dn} t$, lire d, n, t en détachant les lettres.
Là aussi lorsqu'il n'y a pas ambiguité sur le module $k$, on écrit
tout simplement $\mathbf{cn}t$ (resp. $\mathbf{dn}t$) au lieu de
$\mathbf{cn}(t;k)$ (resp. $\mathbf{dn}(t;k)$).

\begin{prop}
On a
\begin{eqnarray}
&&\mathbf{sn}^2 t+\mathbf{cn}^2 t=1,\nonumber\\
&&\mathbf{dn}^2 t+k^2\mathbf{sn}^2 t=1.\nonumber
\end{eqnarray}
\end{prop}
\emph{Démonstration}: En effet, on a
\begin{eqnarray}
\mathbf{cn} t&=&\cos\mathbf{am}t,\nonumber\\
&=&\sqrt{1-\sin^2\mathbf{am}t},\nonumber\\
&=&\sqrt{1-\mathbf{sn}^2 t}.\nonumber
\end{eqnarray}
De même, on a
\begin{eqnarray}
\mathbf{dn} t&=&\sqrt{1-k^2\mathbf{sn}^2 t},\nonumber\\
&=&\sqrt{1-k^2(1-\mathbf{cn}^2 t)},\nonumber\\
&=&\sqrt{1-k^2\mathbf{sn}^2 t}.\nonumber
\end{eqnarray} $\square$

\begin{rem}
Les périodes de la fonction $\mathbf{sn}t$ sont $4K$ et $2iK$ avec
$$K=\int_0^1\frac{ds}{\sqrt{(1-s^2)(1-k^2s^2)}}=F(k),$$
et
$$K'=\int_0^{\frac{1}{k}}\frac{ds}{\sqrt{(s^2-1)(1-k^2s^2)}}
=\int_0^1\frac{ds}{\sqrt{(1-r^2)(1-k'^2r^2)}}=F(k'),$$ où
$k'=\sqrt{1-k^2}$ et $s=\frac{1}{\sqrt{1-k'^2r^2}}$. De même, les
périodes de $\mathbf{cn}t$ sont $4K$ et $2K+2iK'$ et celles de
$\mathbf{dn}t$ sont $2K$ et $4iK'$.
\end{rem}

\begin{prop}
On a
$$\mathbf{sn}(0)=0,\qquad\mathbf{cn}(0)=1,\qquad\mathbf{dn}(0)=1.$$
La fonction $\mathbf{sn}t$ est impaire tandis que les fonctions
$\mathbf{cn}t$ et $\mathbf{dn}t$ sont paires :
$$\mathbf{sn}(-t)=-\mathbf{sn}t,\quad\mathbf{cn}(-t)=\mathbf{cn}t,
\quad\mathbf{dn}(-t)=\mathbf{dn}t.$$
\end{prop}
\emph{Démonstration}: En effet, les trois premières relations sont
évidentes. En ce qui concerne les autres, par définition si
$$t=\int_0^s\frac{ds}{\sqrt{(1-s^2)(1-k^2s^2)}},$$
alors $s=\mathbf{sn}t.$ Dès lors,
$$-t=\int_0^{-s}\frac{ds}{\sqrt{(1-s^2)(1-k^2s^2)}},$$
autrement dit, $\mathbf{sn}(-t)=-s=-\mathbf{sn}t.$ Des relations
$\mathbf{sn}^2t+\mathbf{cn}^2t=1$ et
$\mathbf{dn}^2t+k^2\mathbf{sn}^2t=1,$ on déduit aisément que
$$\mathbf{cn}(-t)=\sqrt{1-\mathbf{sn}^2(-t)}=\sqrt{1-\mathbf{sn}^2t}=\mathbf{cn}t,$$
et
$$\mathbf{dn}(-t)=\sqrt{1-k^2\mathbf{sn}^2(-t)}=\sqrt{1-k^2\mathbf{sn}^2t}=\mathbf{dn}t.$$ $\square$

\begin{prop}
Les dérivées des trois fonctions elliptiques de Jacobi sont
données par
\begin{eqnarray}
\frac{d}{dt}\mathbf{sn}t&=&\mathbf{cn}t.\mathbf{dn}t,\nonumber\\
\frac{d}{dt}\mathbf{cn}t&=&-\mathbf{sn}t.\mathbf{dn}t,\nonumber\\
\frac{d}{dt}\mathbf{dn}t&=&-k^2\mathbf{sn}t.
\mathbf{cn}t.\nonumber
\end{eqnarray}
\end{prop}
\emph{Démonstration}: En effet, par définition si
$$t=\int_0^s\frac{ds}{\sqrt{(1-s^2)(1-k'^2s^2)}},$$
alors $s=\mathbf{sn} t$ et on a
$$\frac{d}{dt}\mathbf{sn}t=\sqrt{(1-\mathbf{sn}^2t)(1-k'^2\mathbf{sn}^2t)}
=\mathbf{cn}t. \mathbf{dn} t.$$ Comme $\mathbf{sn}^2
t+\mathbf{cn}^2 t=1,$ alors
$$\mathbf{sn}t \frac{d}{dt}\mathbf{sn}t+\mathbf{cn}t
\frac{d}{dt}\mathbf{cn}t=0,$$
$$ \mathbf{sn}t \mathbf{cn}t\mathbf{dn}t+\mathbf{cn}t
\frac{d}{dt}\mathbf{cn} t=0,$$
$$\mathbf{sn}t\mathbf{dn}t+\frac{d}{dt}\mathbf{cn}t=0.$$
De même, de la relation $\mathbf{dn}^2t+k^2\mathbf{sn}^2t=1,$ on
déduit que
$$\mathbf{dn}t \frac{d}{dt}\mathbf{dn}t+k^2\mathbf{sn}t\frac{d}{dt}\mathbf{sn}t=0,$$
$$\frac{d}{dt}\mathbf{dn}t+k^2\mathbf{sn}t\mathbf{cn}t=0.$$ $\square$

\begin{prop}
Les fonctions elliptiques de Jacobi vérifient les équations
différentielles :
\begin{eqnarray}
(\frac{d}{dt}\mathbf{sn}t)^2&=&(1-\mathbf{sn}^2t)(1-k^2\mathbf{sn}^2t),\nonumber\\
(\frac{d}{dt}\mathbf{cn}t)^2&=&(1-\mathbf{cn}^2t)(k'^2+k^2\mathbf{cn}^2t),\nonumber\\
(\frac{d}{dt}\mathbf{dn}t)^2&=&(1-\mathbf{dn}^2t)(\mathbf{dn}^2t-k'^2),\nonumber
\end{eqnarray}
où $k'=\sqrt{1-k^2}$.
\end{prop}
\emph{Démonstration}: En effet, la première équation a été obtenue
dans la preuve de la proposition précédente. Concernant les deux
autres équations, on a
\begin{eqnarray}
(\frac{d}{dt}\mathbf{cn}t)^2&=&\mathbf{sn}^2t.\mathbf{dn}^2t,\quad \mbox{(proposition 3.3)}\nonumber\\
&=&(1-\mathbf{cn}^2t)(1-k^2\mathbf{sn}^2t),\quad \mbox{(proposition 3.1)}\nonumber\\
&=&(1-\mathbf{cn}^2t)(1-k^2(1-\mathbf{cn}^2t)),\quad
\mbox{(proposition 3.1)}\nonumber\\
&=&(1-\mathbf{cn}^2t)(k'^2+k^2\mathbf{cn}^2t),\nonumber
\end{eqnarray}
et
\begin{eqnarray}
(\frac{d}{dt}\mathbf{dn}t)^2&=&k^4\mathbf{sn}^2t.\mathbf{cn}^2t,\quad \mbox{(proposition 3.3)}\nonumber\\
&=&k^2(1-\mathbf{dn}^2t)(1-\mathbf{sn}^2t),\quad \mbox{(proposition 3.1)}\nonumber\\
&=&(1-\mathbf{dn}^2t)(k^2-(1-\mathbf{dn}^2t)),\quad
\mbox{(proposition 3.1)}\nonumber\\
&=&(1-\mathbf{dn}^2t)(\mathbf{dn}^2t-k'^2).\nonumber
\end{eqnarray} $\square$

\begin{cor}
Les fonctions elliptiques de Jacobi : $\mathbf{sn}t$,
$\mathbf{cn}t$ et $\mathbf{dn}t$ s'obtiennent par inversion
respectivement des intégrales :
\begin{eqnarray}
t&=&\int_0^w\frac{dw}{\sqrt{(1-w^2)(1-k^2w^2)}},\nonumber\\
t&=&\int_0^w\frac{dw}{\sqrt{(1-w^2)(k'^2+k^2w^2)}},\nonumber\\
t&=&\int_0^w\frac{dw}{\sqrt{(1-w^2)(w^2-k'^2)}}.\nonumber
\end{eqnarray}
où $k'=\sqrt{1-k^2}$.
\end{cor}
\emph{Démonstration}: En posant $w=\mathbf{sn} t$ dans la première
équation différentielle (proposition 3.4), on obtient
$$\frac{dw}{dt}=\sqrt{(1-w^2)(1-k^2w^2)},$$
et il suffit de noter que :
\begin{eqnarray}
w(0)&=&\mathbf{sn}(0)=0,\quad \mbox{(proposition 3.2)}\nonumber\\
\frac{dw}{dt}(0)&=&\mathbf{sn}'(0),\nonumber\\
&=&\mathbf{cn}(0).\mathbf{dn}(0),\quad \mbox{(proposition 3.3)}\nonumber\\
&=&1.\quad \mbox{(proposition 3.2)}\nonumber
\end{eqnarray}
De même, en posant $w=\mathbf{cn} t$ dans la seconde équation
différentielle (proposition précédente), on obtient
$$\frac{dw}{dt}=\sqrt{(1-w^2)(k'^2+k^2w^2)},$$
et il suffit de noter que :
\begin{eqnarray}
w(0)&=&\mathbf{cn}(0)=1,\quad \mbox{(proposition 3.2)}\nonumber\\
\frac{dw}{dt}(0)&=&\mathbf{cn}'(0),\nonumber\\
&=&-\mathbf{sn}(0).\mathbf{dn}(0),\quad \mbox{(proposition 3.3)}\nonumber\\
&=&0.\quad \mbox{(proposition 3.2)}\nonumber
\end{eqnarray}
Et enfin, en posant $w=\mathbf{dn} t$ dans la troisième équation
différentielle (proposition précédente), on obtient
$$\frac{dw}{dt}=\sqrt{(1-w^2)(w^2-k'^2)},$$
et il suffit de noter que :
\begin{eqnarray}
w(0)&=&\mathbf{dn}(0)=1,\quad \mbox{(proposition 3.2)}\nonumber\\
\frac{dw}{dt}(0)&=&\mathbf{dn}'(0),\nonumber\\
&=&-k^2\mathbf{sn}(0).\mathbf{cn}(0),\quad \mbox{(proposition 3.3)}\nonumber\\
&=&0. \quad \mbox{(proposition 3.2)}\nonumber
\end{eqnarray} $\square$

Examinons enfin le cas où $k=0$ et $k=1.$

\begin{prop}
a) Quand $k=0$, on a
$$
\mathbf{am}(t;0)=t,\quad \mathbf{sn}(t;0)=\sin t,\quad
\mathbf{cn}(t;0)=\cos t,\quad \mathbf{dn}(t;0)=1.
$$
b) Lorsque $k=1$, on a
$$
\mathbf{sn}(t;1)=\tanh t,\quad \mathbf{cn}(t;1)=\frac{1}{\cosh
t},\quad \mathbf{cn}(t;1)=\frac{1}{\cosh^2 t}.
$$
\end{prop}
\emph{Démonstration}: En effet, les deux premières relations
s'obtiennent directement en utilisant la définition de ces
intégrales tandis que les autres découlent des relations
$\mathbf{sn}^2t+\mathbf{cn}^2t=1$ et
$\mathbf{dn}^2t+k^2\mathbf{sn}^2t=1$ (proposition 3.1).\\
b) En effet, pour $k=1$ on a
$$t=\int_0^s\frac{ds}{1-s^2}=\frac{1}{2}\ln(\frac{1+s}{1-s})=\arg
\tanh s,\qquad s^2<1,$$ et alors $s=\mathbf{sn}t=\tanh t.$ Pour
les autres relations, on a
$$\mathbf{cn}(t;1)=\sqrt{1-\mathbf{sn}^2(t;1)}=\sqrt{1-\tanh^2 t}=\frac{1}{\cosh
t},$$ et
$$\mathbf{dn}(t;1)=1-\mathbf{sn}^2(t;1)=1-\tanh^2
t=\frac{1}{\cosh^2 t}.$$ $\square$

\begin{prop}
Les fonctions $\mathbf{sn}t$, $\mathbf{cn}t$, $\mathbf{dn}t$
satisfont respectivement aux formules d'addition suivantes :
\begin{eqnarray}
\mathbf{sn}(t+\tau)&=&\frac{\mathbf{sn}t\mathbf{cn}\tau\mathbf{dn}\tau+
\mathbf{sn}\tau\mathbf{cn}t\mathbf{dn}t}{1-k^2\mathbf{sn}^2t\mathbf{sn}^2\tau},\nonumber\\
\mathbf{cn}(t+\tau)&=&\frac{\mathbf{cn}t\mathbf{cn}\tau-
\mathbf{sn}t\mathbf{sn}\tau\mathbf{dn}t\mathbf{dn}\tau}
{1-k^2\mathbf{sn}^2t\mathbf{sn}^2\tau},\nonumber\\
\mathbf{dn}(t+\tau)&=&\frac{\mathbf{dn}t\mathbf{dn}\tau-
k^2\mathbf{sn}t\mathbf{sn}\tau\mathbf{cn}t\mathbf{cn}\tau}
{1-k^2\mathbf{sn}^2t\mathbf{sn}^2\tau}.\nonumber
\end{eqnarray}
\end{prop}
\emph{Démonstration}: Considérons l'équation d'Euler
$$\frac{ds}{\sqrt{P(s)}}+\frac{dr}{\sqrt{P(r)}}=0,$$
où $P(\xi)=(1-\xi^2)(1-k^2\xi^2), 0<k<1$. L'intégrale de cette
équation peut s'écrire sous la forme
\begin{equation}\label{eqn:euler}
t+\tau=C_1,
\end{equation}
où $C_1$ est une constante et
$$t=\int_0^s\frac{ds}{\sqrt{P(s)}},\qquad
\tau=\int_0^r\frac{dr}{\sqrt{P(r)}},$$ avec $s=\mathbf{sn}t$ et
$r=\mathbf{sn}\tau$. Considérons maintenant le système
différentiel
\begin{eqnarray}
\frac{ds}{dz}&=&\sqrt{P(s)},\\
\frac{dr}{dz}&=&\sqrt{P(r)}.\nonumber
\end{eqnarray}
On a
\begin{eqnarray}
\frac{d^2s}{dz^2}&=&s(2k^2s^2-1-k^2),\nonumber\\
\frac{d^2r}{dz^2}&=&r(2k^2r^2-1-k^2),\nonumber
\end{eqnarray}
et
\begin{eqnarray}
r\frac{d^2s}{dz^2}-s\frac{d^2r}{dz^2}&=&2k^2sr(s^2-r^2),\nonumber\\
r^2(\frac{ds}{dz})^2-s^2(\frac{dr}{dz})^2&=&(r^2-s^2)(1-k^2s^2r^2).\nonumber
\end{eqnarray}
Notons que
$$\frac{\frac{d}{dz}(r\frac{ds}{dz}-s\frac{dr}{dz})}
{(r\frac{ds}{dz}+s\frac{dr}{dz})(r\frac{ds}{dz}-s\frac{dr}{dz})}=
\frac{r\frac{d^2s}{dz^2}-s\frac{d^2r}{dz^2}}{r^2(\frac{ds}{dz})^2-s^2(\frac{dr}{dz})^2},$$
donc
$$\frac{\frac{d}{dz}(r\frac{ds}{dz}-s\frac{dr}{dz})}
{r\frac{ds}{dz}-s\frac{dr}{dz}}=\frac{2k^2sr}{k^2s^2r^2-1}(r\frac{ds}{dz}+s\frac{dr}{dz}).$$
En intégrant, on obtient
$$\frac{d}{dz}\ln(r\frac{ds}{dz}-s\frac{dr}{dz})=\frac{d}{dz}\ln(k^2s^2r^2-1),$$
d'où
$$r\frac{ds}{dz}-s\frac{dr}{dz}=C_2(1-k^2s^2r^2),$$
où $C_2$ est liée à $C_1$ par une relation de la forme :
$C_2=f(C_1)$ avec $f$ une fonction à déterminer. En tenant compte
de (3.2), on obtient
\begin{equation}\label{eqn:euler}
r\sqrt{P(s)}+s\sqrt{P(r)}=C_2(1-k^2s^2r^2).
\end{equation}
Or
\begin{eqnarray}
\sqrt{P(s)}&=&\sqrt{(1-s^2)(1-k^2s^2)},\nonumber\\
&=&\sqrt{(1-\mathbf{sn}^2t)(1-k^2\mathbf{sn}^2t)},\nonumber\\
&=&\mathbf{cn}t\mathbf{dn}t,\quad(\mbox{proposition 3.1})\nonumber\\
\sqrt{P(s)}&=&\mathbf{cn}\tau\mathbf{dn}\tau,\nonumber
\end{eqnarray}
donc l'équation (3.3) devient
$$\mathbf{sn}t\mathbf{cn}\tau\mathbf{dn}\tau+
\mathbf{sn}\tau\mathbf{cn}t\mathbf{dn}t=
C_2(1-k^2\mathbf{sn}^2t\mathbf{sn}^2\tau).$$ Rappelons que
$C_2=f(C_1)=f(t+\tau)$ (d'après (3.1)). Dès lors pour $\tau=0$, on
a $f(t)=\mathbf{sn}t$. Donc
$$\frac{\mathbf{sn}t\mathbf{cn}\tau\mathbf{dn}\tau+
\mathbf{sn}\tau\mathbf{cn}t\mathbf{dn}t}{1-k^2\mathbf{sn}^2t\mathbf{sn}^2\tau}
=\mathbf{sn}(t+\tau).$$ Pour les deux autres formules, il suffit
d'utiliser un raisonnement similaire au précédent. $\square$

\begin{rem}
En un certain sens, les fonctions elliptiques de Jacobi
$\mathbf{sn}t$ et $\mathbf{cn}t$ généralisent les fonctions
trigonométriques sinus et cosinus.
\end{rem}

\section{Applications}

\subsection{Le pendule simple}

Le pendule simple est constitué par un point matériel suspendu à
l'extrémité d'un fil (ou une tige théoriquement sans masse)
astreint à se mouvoir sans frottement sur un cercle vertical. On
désigne par $l$ la longueur du fil (i.e., le rayon du cercle), $g$
l'accélération de la pesanteur et $x$ l'angle instantané du fil
avec la verticale. L'équation du mouvement est
\begin{equation}\label{eqn:euler}
\frac{d^2x}{dt^2}+\frac{g}{l}\sin x=0.
\end{equation}
Posons $\theta=\frac{dx}{dt}$, l'équation (4.1) s'écrit
$$\theta d\theta+\frac{g}{l}\sin x dx=0.$$
En intégrant, on obtient
$$\frac{\theta^2}{2}=\frac{g}{l}\cos x+C,$$
où $C$ est une constante. Pour déterminer cette dernière, notons
que lorsque $t=0$, $x=x_0$ (angle initial), alors $\theta=0$ (la
vitesse est nulle), d'où
$$C=-\frac{g}{l}\cos x_0.$$
Par conséquent
\begin{equation}\label{eqn:euler}
\frac{l}{2g}(\frac{dx}{dt})^2=\frac{l}{2g}\theta^2=\cos x-\cos
x_0.
\end{equation}
Nous allons étudier plusieurs cas :\\
a) Considérons le cas d'un mouvement oscillatoire, i.e., le cas où
la masse passe de $x=x_0$ (le plus grand angle atteint par le
pendule; il y correspond une vitesse $\theta=0$) à $x=0$ (vitesse
maximale). Comme $\cos x=1-2\sin^2\frac{x}{2}$, alors l'équation
(4.2) devient
\begin{equation}\label{eqn:euler}
\frac{l}{4g}(\frac{dx}{dt})^2=\sin^2\frac{x_0}{2}-\sin^2\frac{x}{2}.
\end{equation}
Posons
$$\sin\frac{x}{2}=\sin\frac{x_0}{2}\sin \varphi,$$
d'où
$$\frac{1}{2}\cos\frac{x}{2}dx=\sin\frac{x_0}{2}\cos \varphi
d\varphi,$$
$$\frac{1}{2}\sqrt{1-\sin^2\frac{x}{2}}dx=\sin\frac{x_0}{2}\sqrt{1-\sin^2 \varphi}
d\varphi,$$
$$\frac{1}{2}\sqrt{1-\sin^2\frac{x_0}{2}\sin^2\varphi}dx
=\sin\frac{x_0}{2}\sqrt{1-\sin^2 \varphi} d\varphi,$$ et donc
$$dx=\frac{2\sin\frac{x_0}{2}\sqrt{1-\sin^2
\varphi}}{\sqrt{1-\sin^2\frac{x_0}{2}\sin^2\varphi}}d\varphi.$$
Par substitution dans (4.3), on obtient
$$(\frac{d\varphi}{dt})^2=\frac{g}{l}(1-k^2\sin^2\varphi),$$
où $$k=\sin\frac{x_0}{2},$$ est le module et $\frac{x_0}{2}$
l'angle modulaire. Notons que pour $x=0$ on a $\varphi=0$ et dès
lors
$$t=\pm
\sqrt{\frac{l}{g}}\int_0^\varphi\frac{d\varphi}{\sqrt{1-k^2\sin^2\varphi}}.$$
D'après la section 3, on a donc
\begin{eqnarray}
\varphi&=&\pm\mbox{\textbf{am}}\sqrt{\frac{g}{l}}t,\nonumber\\
\sin \varphi&=&\pm \sin \mbox{\textbf{am}}\frac{g}{l}t=\pm
\mbox{\textbf{sn}}\sqrt{\frac{g}{l}}t,\nonumber\
\end{eqnarray}
et par conséquent
$$\sin\frac{x}{2}=\pm \sin\frac{x_0}{2}\mbox{\textbf{sn}}\sqrt{\frac{g}{l}}t.$$
b) Considérons le cas d'un mouvement circulaire. On écrit
l'équation (4.2) sous la forme
\begin{eqnarray}
\frac{l}{2g}(\frac{dx}{dt})^2&=&1-2\sin^2\frac{x}{2}-\cos
x_0,\nonumber\\
&=&(1-\cos x_0)(1-k^2\sin^2\frac{x}{2}),\nonumber
\end{eqnarray}
où $$k^2=\frac{2}{1-\cos x_0},$$ avec $k$ positif et $0<k<1$. En
tenant compte de la condition initiale $x(0)=0$, on obtient
$$dt=\pm \sqrt{\frac{2l}{g(1-\cos
x_0)}}\int_0^\varphi\frac{d\varphi}{\sqrt{1-k^2\sin^2\varphi}},\quad
\varphi=\frac{x}{2}.$$ Donc
$$\varphi=\pm \mbox{\textbf{am}}\sqrt{\frac{g(1-\cos
x_0)}{2l}}t,$$ et
$$x=\pm 2\mbox{\textbf{am}}\sqrt{\frac{g(1-\cos
x_0)}{2l}}t.$$ c) Considérons enfin le cas d'un mouvement
asymptotique. C'est le cas où $x_0=\pm \pi$ et l'équation (4.2)
s'écrit
$$\frac{l}{2g}(\frac{dx}{dt})^2=\cos x+1=2\cos^2\frac{x}{2}.$$
D'où
\begin{eqnarray}
t&=&\pm \frac{1}{2}\sqrt{\frac{l}{g}}\int_0^x\frac{dx}{\cos \frac{x}{2}},\nonumber\\
&=&\pm \sqrt{\frac{l}{g}}\ln \tan
(\frac{x}{\pi}+\frac{\pi}{4}),\nonumber
\end{eqnarray}
et
$$x=4\arctan e^{\pm \sqrt{\frac{g}{l}}t}-\pi.$$
On vérifie que $x\rightarrow \pm \pi$ quand $t\rightarrow \infty$.

\begin{rem}
Pour des petites oscillations, on peut approcher $\sin x$ par $x$
et l'équation (4.1) se ramène à une équation linéaire,
$$\frac{d^2x}{dt^2}+\frac{g}{l}x=0,$$
dont la solution générale est immédiate :
$$x(t)=C_1\cos \sqrt{\frac{g}{l}}t+C_2\sqrt{\frac{l}{g}}\sin
\sqrt{\frac{g}{l}}t,$$ où $C_1=x(0)$ et $C_2=\frac{dx}{dt}(0)$.
Pour des petites oscillations la période du pendule (le temps
nécessité pour une oscillation complète; un aller-retour) est
$2\pi \sqrt{\frac{l}{g}}$. Par contre, dans le cas des
oscillations qui ne sont pas nécessairement petites, la période
vaut d'après ce qui précéde
$4\sqrt{\frac{l}{g}}\int_0^{\frac{\pi}{2}}\frac{dx}{\sqrt{1-k^2\sin
^2x}}$ avec $k=\sin \frac{x_0}{2}$.
\end{rem}

\subsection{Le corps solide d'Euler}

Les équations d'Euler\footnote{On parle aussi de mouvement
d'Euler-Poinsot du solide} du mouvement de rotation d'un solide
autour d'un point fixe, pris comme origine du repère lié au
solide, lorsqu'aucune force extérieure n'est appliquée au système,
peuvent s'écrire sous la forme
\begin{equation}\label{eqn:euler}
\left\{\begin{array}{rl}
&\frac{dm_{1}}{dt}=\left(\lambda_{3}-\lambda _{2}\right)m_{2}m_{3},\\
&\frac{dm_{2}}{dt}=\left( \lambda _{1}-\lambda _{3}\right) m_{1}m_{3},\\
&\frac{dm_{3}}{dt}=\left(\lambda_{2}-\lambda_{1}\right)m_{1}m_{2}.
\end{array}\right.
\end{equation}
où $( m_{1},m_{2},m_{3})$ est le moment angulaire du solide et
$\lambda _{i}\equiv I_{i}^{-1}$, $I_{1},I_{2}$ et $I_{3}$ étant
les moments d'inertie. Ces équations admettent deux intégrales
premières quadratiques :
$$H_1=\frac{1}{2}\left( \lambda _{1}m_{1}^{2}+\lambda
_{2}m_{2}^{2}+\lambda _{3}m_{3}^{2}\right) ,$$
et
$$H_2=\frac{1}{2}\left( m_{1}^{2}+m_{2}^{2}+m_{3}^{2}\right).$$
Nous supposerons que $\lambda_1, \lambda_2, \lambda_3$ sont tous
différents de zero\footnote{c'est-à-dire que le solide n'est pas
réduit à un point et n'est pas non plus concentré sur une
droite.}. Dans ces conditions, $H_1=0$ entraine $m_1=m_2=m_3=0$ et
donc $H_2=0$; le solide est au repos. Nous écartons ce cas trivial
et supposons dorénavant que $H_1\neq 0$ et $H_2\neq 0$. Lorsque
$\lambda_1=\lambda_2=\lambda_3$, les équations (4.4) montrent
évidemment que $m_1$, $m_2$ et $m_3$ sont des constantes.
Supposons par exemple que $\lambda_1=\lambda_2$, les équations
(4.4) s'écrivent alors
\begin{eqnarray}
\frac{dm_{1}}{dt}&=&\left( \lambda _{3}-\lambda _{1}\right)
m_{2}m_{3},\nonumber\\
\frac{dm_{2}}{dt}&=&\left( \lambda _{1}-\lambda _{3}\right) m_{1}m_{3},\nonumber\\
\frac{dm_{3}}{dt}&=&0.\nonumber
\end{eqnarray}
On déduit alors que $m_3=\mbox{constante}\equiv A$ et
\begin{eqnarray}
\frac{dm_{1}}{dt}&=&A\left( \lambda _{3}-\lambda _{1}\right)
m_{2},\nonumber\\
\frac{dm_{2}}{dt}&=&A\left( \lambda _{1}-\lambda _{3}\right)
m_{1}.\nonumber
\end{eqnarray}
Notons que
$$\frac{d}{dt}(m_1+im_2)=iA( \lambda _{1}-\lambda
_{3})(m_1+im_2),$$ on obtient $m_1+im_2=Ce^{iA( \lambda
_{1}-\lambda _{3})t},$ où C est une constante et donc
$$
m_1=C\cos A( \lambda _{1}-\lambda _{3})t,\quad m_2=C\sin A(
\lambda _{1}-\lambda _{3})t
$$
L'intégration des équations d'Euler est délicate dans le cas
général où $\lambda _{1}$, $\lambda _{2}$ et $\lambda _{3}$ sont
tous différents; les solutions s'expriment à l'aide de fonctions
elliptiques. Dans la suite nous supposerons que $\lambda _{1}$,
$\lambda _{2}$ et $\lambda _{3}$ sont tous différents et nous
écartons les autres cas triviaux qui ne posent aucune difficulté
pour la résolution des équations en question. Pour fixer les idées
nous supposerons dans la suite que :
$\lambda_1>\lambda_2>\lambda_3.$ Géométriquement, les équations
\begin{equation}\label{eqn:euler}
\lambda _{1}m_{1}^{2}+\lambda _{2}m_{2}^{2}+\lambda
_{3}m_{3}^{2}=2H_{1},
\end{equation}
et
\begin{equation}\label{eqn:euler}
m_{1}^{2}+m_{2}^{2}+m_{3}^{2}=2H_{2}\equiv r^2,
\end{equation}
représentent respectivement les équations de la surface d'un
ellipsoide de demi-axes : $\sqrt{\frac{2H_1}{\lambda_1}}$ (demi
grand axe),$\sqrt{\frac{2H_1}{\lambda_2}}$(demi axe
moyen),$\sqrt{\frac{2H_1}{\lambda_3}}$(demi petit axe), et d'une
sphère de rayon r. Donc le mouvement du solide s'effectue sur
l'intersection d'un ellipsoide avec une sphère. Cette intersection
a un sens car en comparant (4.5) à (4.6), on voit que
$\frac{2H_1}{\lambda_1}<r^2<\frac{2H_1}{\lambda_3},$ ce qui
signifie géométriquement que le rayon de la sphère (4.6) est
compris entre le plus petit et le plus grand des demi-axes de
l'ellipsoïde (4.5). Pour étudier l'allure des courbes
d'intersection de l'éllipsoïde (4.5) avec la sphère (4.6), fixons
$H_1>0$ et faisons varier le rayon r. Comme
$\lambda_1>\lambda_2>\lambda_3,$ les demi-axes de l'ellipsoïde
seront
$\frac{2H_1}{\lambda_1}>\frac{2H_1}{\lambda_2}>\frac{2H_1}{\lambda_3}.$
Si le rayon r de la sphère est inférieur au demi petit axe
$\frac{2H_1}{\lambda_3}$ ou supérieur au demi grand axe
$\frac{2H_1}{\lambda_1}$, alors l'intersection en question est
vide ( et aucum mouvement réel ne correspond à ces valeurs de
$H_1$ et r). Lorsque le rayon r est égal à
$\frac{2H_1}{\lambda_3},$ alors l'intersection est composée de
deux points. Lorsque le rayon r augmente
($\frac{2H_1}{\lambda_3}<r<\frac{2H_1}{\lambda_2}$), on obtient
deux courbes autour des extrémités du demi petit axe. De même si
$r=\frac{2H_1}{\lambda_1},$ on obtient les deux extrémités du demi
grand axe et si r est légérement inférieur à
$\frac{2H_1}{\lambda_1},$ on obtient deux courbes fermées au
voisinage de ces extrémités. Enfin, si $r=\frac{2H_1}{\lambda_2}$
alors l'intersection en question est constituée de deux cercles.

\begin{prop}
Les équations différentielles (4.4) d'Euler, s'intégrent au moyen
de fonctions elliptiques de Jacobi.
\end{prop}
\emph{Démonstration}: A partir des intégrales premières (4.5) et
(4.6), on exprime $m_1$ et $m_3$ en fonction de $m_2$. On
introduit ensuite ces expressions dans la seconde \'{e}quation du
système (4.4) pour obtenir une équation différentielle en $m_2$ et
$\frac{dm_2}{dt}$ seulement. De manière plus détaillée, on tire
aisément de (4.5) et (4.6) les relations suivantes
\begin{eqnarray}
m_{1}^2&=&\frac{2H_1-r^2\lambda _{3}- \left( \lambda _{2}-\lambda
_{3}\right) m_{2}^{2}}{\lambda _{1}- \lambda _{3}},\\
m_{3}^2&=&\frac{r^2\lambda _{1}-2H_1-\left( \lambda _{1}-\lambda
_{2}\right) m_{2}^{2}}{\lambda _{1}-\lambda _{3}}.
\end{eqnarray}
En substituant ces expressions dans la seconde équation du système
(4.4), on obtient
$$\frac{dm_{2}}{dt}=\sqrt{(2H_1-r^2\lambda _{3}-
\left( \lambda _{2}-\lambda _{3}\right) m_{2}^{2})(r^2\lambda
_{1}-2H_1-\left( \lambda _{1}-\lambda _{2}\right) m_{2}^{2}) }.$$
En intégrant cette équation, on obtient une fonction $t(m_2)$ sous
forme d'une intégrale elliptique. Pour réduire celle-ci à la forme
standard, on peut supposer que $r^2>\frac{2H_1}{\lambda_2}$
(sinon, il suffit d'intervertir les indices 1 et 3 dans toutes les
formules précédentes). On réecrit l'équation précédente, sous la
forme
$$\frac{dm_{2}}{\sqrt{(2H_1-r^2\lambda _{3})(r^2\lambda
_{1}-2H_1)}dt}= \sqrt{(1- \frac{\lambda _{2}-\lambda
_{3}}{2H_1-r^2\lambda _{3}}m_{2}^{2})(1-\frac{\lambda _{1}-\lambda
_{2}}{r^2\lambda _{1}-2H_1} m_{2}^{2})}.$$ En posant
\begin{eqnarray}
\tau&=&t\sqrt{(\lambda _{2}-\lambda _{3})(r^2\lambda
_{1}-2H_1)},\nonumber\\
s&=&m_2\sqrt{\frac{\lambda _{2}-\lambda _{3}}{2H_1-r^2\lambda
_{3}}},\nonumber
\end{eqnarray}
on obtient
$$\frac{ds}{d\tau}=
\sqrt{(1- s^{2})(1-\frac{(\lambda _{1}-\lambda
_{2})(2H_1-r^2\lambda _{3})}{(\lambda_2-\lambda_3)(r^2\lambda
_{1}-2H_1)} s^{2})},$$ ce qui suggère de choisir comme module des
fonctions elliptiques
$$k^2=\frac{(\lambda _{1}-\lambda
_{2})(2H_1-r^2\lambda _{3})}{(\lambda_2-\lambda_3)(r^2\lambda
_{1}-2H_1)}.$$ Les inégalités $\lambda_1>\lambda_2>\lambda_3$,
$\frac{2H_1}{\lambda_1}<r^2<\frac{2H_1}{\lambda_3}$ et
$r^2>\frac{2H_1}{\lambda_2}$ montrent qu'effectivement $0<k^2<1.$
On obtient donc
$$\frac{ds}{d\tau}=
\sqrt{(1- s^{2})(1-k^2 s^{2})}.$$ Cette équation admet la
solution\footnote{on convient de choisir l'origine des temps telle
que $m_2=0$ pour $t=0.$}
$$\tau=\int_0^s\frac{ds}{\sqrt{(1-s^2)(1-k^2s^2)}}.$$
La fonction inverse $s(\tau)$ constitue l'une des fonctions
elliptiques de Jacobi : $s=\mathbf{sn} \tau,$ qui détermine
également $m_2$ en fonction du temps, i.e.,
$$
m_2=\sqrt{\frac{2H_1-r^2\lambda_{3}}{\lambda_2-\lambda_3}}\cdot\mathbf{sn}\tau.
$$
D'après les égalités (4.7) et (4.8), on sait que les fonctions
$m_1$ et $m_3$ s'expriment algébriquement à l'aide de $m_2,$ donc
$$m_1=\sqrt{\frac{2H_1-r^2\lambda_{3}}{\lambda_1-\lambda_3}}
\cdot\sqrt{1-\mathbf{sn}^2\tau},$$ et
$$m_3=\sqrt{\frac{r^2\lambda_{1}-2H_1}{\lambda_1-\lambda_3}}\cdot
\sqrt{1-k^2\mathbf{sn}^2\tau}.$$ Compte tenu de la définition des
deux autres fonctions elliptiques (voir section 3)
$$\mathbf{cn}\tau=\sqrt{1-\mathbf{sn}^2\tau},\qquad \mathbf{dn}\tau=
\sqrt{1-k^2\mathbf{sn}^2\tau},$$ et du fait que
$\tau=t\sqrt{(\lambda _{2}-\lambda _{3})(r^2\lambda _{1}-2H_1)}$,
on obtient finalement les formules suivantes :
\begin{equation}\label{eqn:euler}
\left\{\begin{array}{rl}
m_1=\sqrt{\frac{2H_1-r^2\lambda_{3}}{\lambda_1-\lambda_3}}&\mathbf{cn}(t\sqrt{(\lambda
_{2}-\lambda _{3})(r^2\lambda
_{1}-2H_1)}),\\
m_2=\sqrt{\frac{2H_1-r^2\lambda_{3}}{\lambda_2-\lambda_3}}&\mathbf{sn}(t\sqrt{(\lambda
_{2}-\lambda _{3})(r^2\lambda
_{1}-2H_1)}),\\
m_3=\sqrt{\frac{r^2\lambda_{1}-2H_1}{\lambda_1-\lambda_3}}&\mathbf{dn}(t\sqrt{(\lambda
_{2}-\lambda _{3})(r^2\lambda _{1}-2H_1)}).
\end{array}\right.
\end{equation}
Autrement dit, l'intégration des équations d'Euler s'effectue au
moyen de fonctions elliptiques. $\square$

\begin{rem}
Notons que pour $\lambda_1=\lambda_2$, on a $k^2=0$. Dans ce cas,
les fonctions elliptiques $\mathbf{sn}\tau, \mathbf{cn}\tau,
\mathbf{dn}\tau$ se réduisent respectivement aux fonctions
$\sin\tau, \cos\tau, 1$. Dès lors de (4.9), on tire aisément que
$$
\left\{\begin{array}{rl}
m_1=\sqrt{\frac{2H_1-r^2\lambda_{3}}{\lambda_1-\lambda_3}}&
\cos\sqrt{(\lambda_1-\lambda_3)(r^2\lambda_{1}-2H_1)}t,\\
m_2=\sqrt{\frac{2H_1-r^2\lambda_{3}}{\lambda_1-\lambda_3}}&
\sin\sqrt{(\lambda_1-\lambda_3)(r^2\lambda_{1}-2H_1)}t,\\
m_3=\sqrt{\frac{r^2\lambda_{1}-2H_1}{\lambda_1-\lambda_3}}&.
\end{array}\right.
$$
On retrouve les solutions établis précédemment avec
$A=\sqrt{\frac{r^2\lambda_{1}-2H_1}{\lambda_1-\lambda_3}}$ et $
C=\sqrt{\frac{2H_1-r^2\lambda_{3}}{\lambda_1-\lambda_3}}.$
\end{rem}

\subsection{Une famille de syst\`{e}mes int\'{e}grables}

On consid\`{e}re un syst\`{e}me diff\'{e}rentiel non-lin\'{e}aire
sur $\mathbb{R}^{4}$ d\'{e}fini par le hamiltonien
\begin{equation}\label{eqn:euler}
H=\frac{1}{2}[
x_{1}^{2}+x_{2}^{2}+a(y_1^2+y_2^2)+b(y_1^2+y_2^2)^2+c(y_1^2+y_2^2)^3]
,
\end{equation}
où $a$, $b$, $c$, sont des constantes. Dans ce cas, le syst\`{e}me
dynamique hamiltonien associ\'{e} \`{a} $H$ s'\'{e}crit
\begin{eqnarray}
\frac{dy_{1}}{dt}&=&x_{1},\nonumber\\
\frac{dy_{2}}{dt}&=&x_{2},\\
\frac{dx_{1}}{dt}&=&-[a+2b(y_1^2+y_2^2)+3c(y_1^2+y_2^2)^2]y_1,\nonumber\\
\frac{dx_{2}}{dt}&=&-[a+2b(y_1^2+y_2^2)+3c(y_1^2+y_2^2)^2]y_2.\nonumber
\end{eqnarray}
Ces \'{e}quations donnent un champ de vecteurs sur
$\mathbb{R}^{4}$.

\begin{prop}
Le système (4.11) admet une intégrale première quadratique qui
détermine avec $H$(4.10) un système intégrable au sens de
Liouville. En outre, la linéarisation c'est-à-dire la description
des niveaux communs des intégrales et les flots dont ils sont
pourvus s'effectue sur une courbe elliptique.
\end{prop}
\emph{Démonstration}: Dans notre cas, l'existence d'une seconde
intégrale première indépendante et en involution avec $H_{1}\equiv
H,$ suffit pour que le système soit intégrable au sens de
Liouville (voir appendice). Le système différentiel (4.11)
implique
$$\frac{d^2y_1}{dt^2}+[a+2b(y_1^2+y_2^2)+3c(y_1^2+y_2^2)^2]y_1=0,$$
$$\frac{d^2y_2}{dt^2}+[a+2b(y_1^2+y_2^2)+3c(y_1^2+y_2^2)^2]y_2=0,$$
d'où $$y_2\frac{d^2y_1}{dt^2}-y_1\frac{d^2y_2}{dt^2}=0.$$ Dès lors
$$\frac{d}{dt}(y_2\frac{dy_1}{dt}-y_1\frac{dy_2}{dt})=0,$$
de sorte que la fonction (le moment)
$$H_2=x_{1}y_{2}-x_{2}y_{1},$$ est une intégrale première.
Les fonctions $H_1$ et $H_2$ sont en involution $$\{H_1,H_2\}
=\sum_{i=1}^{2}(\frac{\partial H_1}{\partial x_k}\frac{\partial
H_2 }{\partial y_k}-\frac{\partial H_1}{\partial
y_k}\frac{\partial H_2}{\partial x_k})=0.$$ Donc cette seconde
intégrale première, détermine avec $H_1$ un système intégrable au
sens de Liouville. Soit $$\{x\equiv (y_1,y_2,x_1,x_2)\in
\mathbb{R}^4:H_1(x)=c_1,H_2(x)=c_2\},$$ la surface invariante où
$(c_1,c_2)$ n'est pas une valeur critique. En substituant
$$y_1=r\cos\theta, \qquad y_2=r\sin\theta,$$ dans les équations
\begin{eqnarray}
H_1&=&\frac{1}{2}[
x_{1}^{2}+x_{2}^{2}+a(y_1^2+y_2^2)+b(y_1^2+y_2^2)^2+c(y_1^2+y_2^2)^3]=c_1, \nonumber\\
H_2&=&x_{1}y_{2}-x_{2}y_{1}=c_2,\nonumber
\end{eqnarray}
on obtient
\begin{eqnarray}
(\frac{dr}{dt})^2+r^{2}(\frac{d\theta}{dt})^2
+ar^2+br^4+cr^6&=&2c_{1},\nonumber\\
r^{2}\frac{d\theta }{dt}&=&-c_{2}.\nonumber
\end{eqnarray}
D'où
$$(r\frac{dr}{dt})^{2}+ar^4+br^6+cr^8-2c_{1}r^{2}+c_{2}^{2}=0,$$ et par conséquent
$$w^{2}+az^2+bz^3+cz^4-2c_{1}z+c_{2}^{2}=0,$$ où $w=r\frac{dr}{dt},\quad z=r^{2}$. La courbe algébrique
$$\mathcal{C}=\overline{\{(w,z):w^{2}+az^2+bz^3+cz^4-2c_{1}z+c_{2}^{2}=0\}},$$
est une courbe elliptique. On a une seule différentielle
holomorphe $$
\omega=\frac{dz}{\sqrt{az^2+bz^3+cz^4-2c_{1}z+c_{2}^{2}}},$$ et la
linéarisation s'effectue donc sur cette courbe elliptique;
autrement dit les équations différentielles (4.11) s'intégrent au
moyen de fonctions elliptiques. $\square$

\subsection{Équations aux d\'{e}riv\'{e}es
partielles coupl\'{e}es non-lin\'{e}aires de Schr\"{o}dinger}

Considérons les équations couplées non linéaires de Schrödinger :
\begin{eqnarray} i\frac{\partial u}{\partial s}=\frac{\partial ^{2}u}{\partial
t^{2}}+(|u|^{2}+|v|^{2})u,\\
i\frac{\partial v}{\partial s}=\frac{\partial ^{2}v}{\partial
t^{2}}+(|u|^{2}+|v| ^{2})v.\nonumber
\end{eqnarray}
Les fonctions $u(s,t)$ et $v(s,t)$ dépendent des variables $s$ et
$t$. On cherche les solutions de $(4.12)$ sous la forme
\begin{eqnarray}
u(s,t)&=&\zeta (t)\exp (ias) ,\nonumber\\
v(s,t)& =&\eta (t)\exp (ias),\nonumber
\end{eqnarray}
où $\zeta(t)$ et $\eta (t)$ sont deux fonctions réelles et $a$\
une constante arbitraire, ce qui implique comme conséquence qu'on
aura
$$\frac{d^{2}\zeta }{dt^{2}}+(a+\zeta ^{2}+\eta ^{2})
\zeta =0,$$ $$\frac{d^{2}\eta }{dt^{2}}+(a+\zeta ^{2}+\eta ^{2})
\eta =0.$$ En posant $$y_{1}=\zeta ,\quad y_{2}=\eta ,\quad
x_{1}=\frac{d\zeta }{dt},\quad x_{2}=\frac{d\eta }{dt},$$ on
obtient
\begin{eqnarray}
\frac{dy_{1}}{dt}&=&x_{1},\nonumber\\
\frac{dy_{2}}{dt}&=&x_{2},\\
\frac{dx_{1}}{dt}&=&-(a+y_{1}^{2}+y_{2}^{2}) y_{1},\nonumber\\
\frac{dx_{2}}{dt}&=&-( a+y_{1}^{2}+y_{2}^{2})y_{2}.\nonumber
\end{eqnarray}
Ces équations donnent un champ de vecteurs sur $\mathbb{R}^4$ et
il apparait ainsi que la résolution du système se trouve ramenée à
la recherche des solutions d'un système dynamique hamiltonien de
la forme $(4.10)$ avec $b=\frac{1}{2}$ et $c=0$. On peut donc
utiliser le résultat obtenu dans la proposition 4.2. Cependant,
nous allons procéder différemment et montrer que le système en
question possède une seconde intégrale première quartique et la
résolution du problème s'effectue aussi en terme de fonctions
elliptiques. Nous allons utiliser la théorie des déformations
isospectrales c'est-à-dire laissant invariant le spectre
d'opérateurs linéaires contenant une indéterminée rationnelle. Les
équations à étudier peuvent être exprimées en termes de relations
de commutation (paire de Lax). Plus précisement, on a le résultat
suivant :
\begin{prop}
Le système différentiel (4.13) admet une paire de Lax de sorte que
la fonction
$$H_2=\frac{a}{2}( x_1^2+x_2^2+a( y_1^2+y_2^2) +\frac{1}{2}(
y_1^2+y_2^2) ^2) +\frac{1}{4}( x_1y_2-x_2y_1) ^2,$$ est une
intégrale première quartique et la linéarisation s'effectue à
l'aide de fonctions elliptiques.
\end{prop}
Démonstration: Considérons la forme de Lax $$\frac{d}{dt}A_{h}=[
B_{h},A_{h}] \equiv B_{h}A_{h}-A_{h}B_{h},$$ où $A_h$ et $B_h$
sont des matrices dépendant d'un paramètre complexe $h$ (paramètre
spectrale). Les coefficients du polynôme caratéristique
$\det(A_{h}-\lambda I),$ ne dépendent pas du temps et ce sont des
intégrales premières en involution. En outre, d'après la méthode
de linéarisation de van Moerbeke-Mumford $[18,13] $ le flot se
linéarise sur un tore algébrique complexe. Celui-ci étant
engendr\'{e} par le réseau définit par la matrice des périodes de
la courbe spectrale d'équation affine
\begin{equation}\label{eqn:euler}
P(h,\lambda) \equiv \det(A_{h}-\lambda I)=0,
\end{equation}
et cette équation décrit une déformation isospectrale. Dans le cas
de notre système, on choisit
$$
A_h= \left(
\begin{array}{cc}
U_h & V_h \\
W_h & -U_h
\end{array}
\right), \qquad
B_h=\left(\begin{array}{cc} 0&1\\
R_h&0\end{array}\right),$$ avec
\begin{eqnarray}
U_{h}&=&\frac{1}{2}(\frac{x_{1}y_{1}+x_{2}y_{2}}{a+h}),\nonumber\\
V_{h}&=&-1-\frac{y_{1}^{2}+y_{2}^{2}}{2(a+h)},\nonumber\\
W_{h}&=&\frac{1}{2}(\frac{x_{1}^{2}+x_{2}^{2}}{a+h})-h+\frac{1}{2}(y_{1}^{2}+y_{}^{2}),\nonumber\\
R_{h}&=&h-y_{1}^{2}-y_{2}^{2}.\nonumber
\end{eqnarray}
Explicitement, l'équation $(4.14)$ fournit
\begin{equation}\label{eqn:euler}
w^{2}=h^{3}+2ah^{2}+(a^{2}-H_{1})h-H_{2},\end{equation} où
\begin{eqnarray}
w&=&\lambda (h+a),\nonumber\\
H_{1}&=&\frac{1}{2}(x_{1}^{2}+x_{2}^{2})+\frac{a}{2}(
y_{1}^{2}+y_{2}^{2})+\frac{1}{4}(
y_{1}^{2}+y_{2}^{2})^{2},\nonumber\\
H_{2}&=&\frac{a}{2}(x_{1}^{2}+x_{2}^{2}+a(
y_{1}^{2}+y_{2}^{2})+\frac{1}{2}(y_{1}^{2}+y_{2}^{2})
^{2})+\frac{1}{4}(x_{1}y_{2}-x_{2}y_{1})^{2},\nonumber\\
 &=&aH_{1}+\frac{1}{4}(x_{1}y_{2}-x_{2}y_{1})^{2}.\nonumber
\end{eqnarray}
Les deux intégrales premières $H_{1}$ et $H_{2}$ sont évidemment
en involution et le système en question est intégrable au sens de
Liouville (voir appendice). Le flot se linéarise sur la courbe
elliptique d'équation affine (4.15). Autrement dit, la
linéarisation s'effectue à l'aide de fonctions elliptiques.
$\square$

\subsection{Le champ de Yang-Mills avec groupe de jauge $SU(2)$}

Soit $F_{kl}$ le champ de Yang-Mills dans l'alg\`{e}bre de Lie
$T_{e}SU(2)$ du groupe $SU(2).$ C'est une expression locale du
champ de Jauge ou connexion $A_{k}$ d\'{e}finissant la
d\'{e}riv\'{e}e covariante de $F_{kl}$ \`{a} l'aide de
l'expression:
$$\triangledown _{k}F_{kl}=\frac{\partial F_{kl}}{\partial \tau _{k}}+
\left[ A_{k},F_{kl}\right] =0,\qquad F_{kl},A_{k}\in
T_{e}SU(2),\quad 1\leq k,l\leq 4,$$ ans laquelle $\left[
A_{k},F_{kl}\right] $ est le crochet des deux champs dans
l'alg\`{e}bre de Lie du groupe de Lie $SU(2)$ et
$$F_{kl}=\frac{\partial A_{l}}{\partial \tau _{k}}-
\frac{\partial A_{k}}{\partial \tau _{l}}+\left[
A_{k},A_{l}\right].$$ Dans le cas qui nous int\'{e}resse, on a
\begin{eqnarray}
\frac{\partial A_{l}}{\partial \tau _{k}}&=&0,\quad \left( k\neq
1\right),\nonumber\\
A_{1}&=&A_{2}=0,\nonumber\\
A_{3}&=&n_{1}U_{1}\in su\left(
2\right),\nonumber\\
A_{4}&=&n_{2}U_{2}\in su\left( 2\right),\nonumber
\end{eqnarray}
o\`{u} $$n_{1}=[n_{2},[n_{1},n_{2}]],\qquad
n_{2}=[n_{1},[n_{2},n_{1}]],$$ engendre $su\left( 2\right) $ et le
syst\`{e}me de Yang-Mills devient
\begin{eqnarray}
\frac{\partial ^{2}U_{1}}{\partial
t^{2}}+U_{1}U_{2}^{2}&=&0,\nonumber\\
\frac{\partial ^{2}U_{2}}{\partial
t^{2}}+U_{2}U_{1}^{2}&=&0,\nonumber
\end{eqnarray}
avec $t=\tau _{1}.$ En posant $U_{1}=q_{1},$ $U_{2}=q_{2},$
$\frac{\partial U_{1}}{\partial t} =p_{1},\frac{\partial
U_{2}}{\partial t}=p_{2},$ les \'{e}quations de Yang-Mills
s'\'{e}crivent sous la forme d'un champ de vecteurs hamiltonien
avec $H=\frac{1}{2}\left(
p_{1}^{2}+p_{2}^{2}+q_{1}^{2}q_{2}^{2}\right)$ l'hamiltonien.
Celui-ci joue un r\^{o}le important en th\'{e}orie des champs. En
utilisant la transformation symplectique
\begin{eqnarray}
p_{1}&=&\frac{\sqrt{2}}{2} \left( x_{1}+x_{2}\right) ,\nonumber\\
p_{2}&=&\frac{\sqrt{2}}{2} \left( x_{1}-x_{2}\right),\nonumber\\
q_{1}&=&\frac{1}{2} \left( \root{4}\of{2}\right) ^{3} \left( y_{1}+iy_{2}\right),\nonumber\\
q_{2}&=&\frac{1}{2} \left( \root{4}\of{2}\right) ^{3} \left(
y_{1}-iy_{2}\right),\nonumber
\end{eqnarray}
on r\'{e}ecrit le hamiltonien ci-dessus sous la forme
$$H=\frac{1}{2}\left( x_{1}^{2}+x_{2}^{2}\right)
+\frac{1}{4}\left( y_{1}^{2}+y_{2}^{2}\right) ^{2},$$ lequel
coincide \'{e}videmment avec $\left(4.10\right) $ pour $a=c=0$,
$b=\frac{1}{2}$ ou ce qui revient au même avec le système
différentiel correspondant $(4.13)$ avec $a=0$. Donc ici aussi
l'intégration du problème en question s'effectue en termes de
fonctions elliptiques.

\subsection{Appendice}

Les équations canoniques de Hamilton s'écrivent sous la forme
\begin{eqnarray}
\frac{dy_{1}}{dt}&=&\frac{\partial H}{\partial x_{1}},\nonumber\\
&\vdots& \nonumber\\
\frac{dy_{n}}{dt}&=&\frac{\partial H}{\partial x_{n}},\\
\frac{dx_{1}}{dt}&=&-\frac{\partial H} {\partial
y_{1}},\nonumber\\
&\vdots& \nonumber\\
\frac{dx_{n}}{dt}&=&-\frac{\partial H}{\partial y_{n}},\nonumber
\end{eqnarray}
où $(x_1,...,x_n)\in\mathbb{R}^n$ et
$(y_1,...,y_n)\in\mathbb{R}^n$, sont des coordonnées dans l'espace
de phase $\mathbb{R}^{2n}$. Ce sont $2n$ équations différentielles
du premier ordre qui sont connues lorsqu'on connait la fonction
$H$ appelée hamiltonien du système. Le système (4.16) est
intégrable au sens de Liouville lorsqu'il possède $n$ intégrales
premières $H_{1}\equiv H,$ $H_{2},\ldots ,H_{n}$ en involution
(c'est-à-dire que les crochets de Poisson
$$\{H_i,H_j\}=\sum_{k=1}^{n}
(\frac{\partial H_i}{\partial x_k}\frac{\partial H_j}{\partial
y_k} -\frac{\partial H_i}{\partial y_k}\frac{\partial
H_j}{\partial x_k}),\quad1\leq i,j\leq n.$$ s'annulent deux à
deux) et qu'en outre les gradients $gradH_{i}$ sont linéairements
indépendants. Pour des constantes génériques $c=(c_1,\ldots ,c_n)
,$ l'ensemble de niveau commun aux intégrales $H_{1},\ldots
,H_{n}$:$$M_c=\{x\equiv (y_1,...,y_n,x_1,...,x_n)\in
\mathbb{R}^{2n}:H_1(x)=c_1,\ldots,H_n(x)=c_n\},$$ forme une
variété de dimension $n$. D'après le théorème d'Arnold-Liouville
$[6, 14]$, si la variété $M_c$ est compacte et connexe, alors elle
est difféomorphe à un tore de dimension $n$:
$$T^n={\mathbb{R}^n}/{\mathbb{Z}^n}=\{(\varphi_1,\ldots,\varphi_n)mod2\pi\},$$
sur lequel le problème se linéarise. En outre, on démontre
l'existence d'une transformation canonique vers de nouvelles
coordonnées, dites variables action-angle, les coordonnées action
étant des constantes du mouvement et les coordonnées angle des
fonctions linéaires dans le temps. Le point $$(y_1(t),\ldots
,y_n(t),x_1(t),\ldots ,x_n(t)),$$ représentant la solution du
système (4.16) a un mouvement quasi-périodique, c'est-à-dire en
coordonnées angulaires $(\varphi_1,\ldots,\varphi_n)$, on a
$$\frac{d\varphi}{dt}=\omega,
\quad \omega=\omega(c)=constante.$$ Donc, en principe, la
transformation canonique fournit les positions et les moments en
fonction du temps et le problème est résolu. La résolution
explicite de plusieurs équations (notamment en mécanique) des
flots hamiltoniens associés aux fonctions $H_1,...,H_n$ se fait à
l'aide d'intégrales elliptiques. Autrement dit pour une
compactification appropriée, on a $\overline{M}_c \simeq
Jac(\mathcal{C)}\simeq\mathcal{C}$ où $\mathcal{C}$ est une courbe
elliptique à déteminer et le système en question s'intégre en
terme de fonctions elliptiques.

\end{document}